\documentclass[times,sort&compress,3p]{elsarticle}
\usepackage{amsmath, amssymb, enumerate, amsthm}

\usepackage[english]{babel}
\usepackage[latin1]{inputenc}
\usepackage{booktabs}
\usepackage{graphicx}
\usepackage{graphics}
\usepackage{caption}
\usepackage{amsmath}
\usepackage{amssymb}
\usepackage{amsthm}
\usepackage{hyperref}
\usepackage{array}
\usepackage{color}
 \usepackage[]{tabularx}

\theoremstyle{plain}% default
\newtheorem{theo}{Theorem}
\newtheorem{lem}[theo]{Lemma}

\newtheorem{corollary}[theo]{Corollary}

\theoremstyle{definition}

\theoremstyle{remark}

\newcommand{\E}{{\rm E}}
\newcommand{\p}{\Pr}
\newcommand{\R}{\mathbb{R}}

\def\NN{\mathcal{N}}

\begin{document}

\begin{frontmatter}

\title{Central limit theorem and  bootstrap procedure for Wasserstein's variations \\
with  an application to structural relationships between distributions}

\author{Eustasio del Barrio$^{a}$, Paula Gordaliza$^{a}$, H\'el\`ene Lescornel$^{b}$ and Jean-Michel 
Loubes$^{b}$\footnote{Corresponding author: loubes@math.univ-toulouse.fr}\\
$^{a}$\textit{IMUVA, Universidad de Valladolid} and $^{b}$\textit{Institut de math\'ematiques de Toulouse}}

\begin{abstract}
Wasserstein barycenters and variance-like  {criteria based on the} Wasserstein distance are used in many problems to analyze the homogeneity of collections of distributions and structural relationships between the observations. 
We propose the estimation of the quantiles of the empirical process of Wasserstein's variation  using a
bootstrap procedure.  {We then} use these results for statistical inference on a distribution registration model 
for general deformation functions.  The tests are based on the variance of the distributions with respect to 
their Wasserstein's barycenters for which we prove central limit theorems, including bootstrap versions.
\end{abstract}


\end{frontmatter}

\section{Introduction\label{s:intro}}

Analyzing the variability of large data sets is a difficult task when  {the inner geometry of} the information conveyed by the observations  {is} far from  {being} Euclidean. Indeed, deformations  on the data such as 
 {location-scale transformations} %a model is not a transformation!
or more general warping procedures  {preclude} the use of  {common} %to avoid using usual...
 {statistical} methods. Looking for a way to measure structural relationships  {within} %not between 
data is of high importance.  {Such issues} arise when considering the estimation of probability measures observed with deformations;  
 {it is common, e.g.,} % in biology, genetic studies are not regarded as biology
when considering gene expression. 

 {Over the last decade, there} has been a large amount of work  {dealing} with registrations issues. We refer,  {e.g.,} to \cite{allassonniere2007towards, JASA, Ramsay-Silverman-05} and references therein. However, when dealing with the registration of warped distributions, the literature is scarce. We mention 
here the method provided for biological computational issues known as quantile normalization in \cite{Bold, GALLON-2011-593476} 
and references therein. Recently, using optimal transport methodologies, comparisons of 
distributions have been studied using a notion of Fr\'echet mean for distributions as in 
\cite{agueh2010barycenters} or a notion of depth as in \cite{2014arXiv1412.8434C}.

 {As a} natural frame for applications  {of a deformation model, consider} $J$ independent  {random samples of size $n$, where for each $j \in \{ 1, \ldots, J \}$, the real-valued random variable $X_j$ has distribution $\mu_j$ and, for each $i \in \{ 1, \ldots, n\}$,  the $i$th observation of $X_j$ is such that} 
$$ 
X_{i,j}=  g_j (\varepsilon_{i,j}),
$$  
where the $\varepsilon_{i,j}$s are iid random variables with unknown distribution $\mu$.  {Assume that the} functions $g_1, \ldots, g_J$ belong to a class $\mathcal{G}$ of deformation functions, which model how the distributions $\mu_1, \ldots, \mu_J$ are warped one to another. 

This model is the natural extension of the functional deformation models studied in the statistical literature for which estimation procedures are provided in \cite{Gamboa-Loubes-Maza-07}  {and} testing issues are tackled in \cite{MR3323102}. Note that at the era of parallelized inference where a large amount of data is processed in the same way but at different locations or by different computers, this framework appears also natural since this parallelization  may lead to small changes with respect to the law of the observations that should be  {eliminated}. 

In the  {framework} of warped distributions, a  {central} goal is the estimation of the warping functions, possibly as a first step towards registration or alignment of the (estimated) distributions. Of course, without some constraints on the class $\mathcal{G}$, the deformation model is meaningless. We can, for instance, obtain any distribution on $\mathbb{R}^d$ as a warped version of a fixed probability having a density if we take the optimal transportation map as the warping function; see 
\cite{villani2009optimal}. One has to consider smaller classes of deformation functions to perform a reasonable registration.

In  {cases where $\mathcal{G}$ is a parametric class}, estimation of the warping functions is studied in \cite{agullo2015parametric}. However,  estimation/registration procedures may lead to inconsistent conclusions if the chosen deformation class $\mathcal{G}$ is too small. It is, therefore, important to be able to assess  {the} fit to the deformation model given by a particular choice of $\mathcal{G}$.  {This} is the main goal of this paper. We note that within this framework, statistical inference on deformation models for distributions has been studied first in \cite{Freitag2005123}. Here we provide a different approach which allows to deal with more general deformation classes.

The  {pioneering} works \cite{MR1704844, MR1625620} study the existence of relationships between distributions $F$ and $G$ by using a discrepancy measure $\Delta(F,G)$ between  {them which is} built using the Wasserstein distance. The authors consider the assumption $\mathcal{H}_0: \Delta(F,G) > \Delta_0$ versus $\mathcal{H}_a: \Delta(F,G) \leq \Delta_0$ for a chosen threshold $\Delta_0$. Thus  when the  {null hypothesis} is rejected, there is statistical evidence that the two distributions are similar with respect to the chosen criterion. In this  {same vein}, we define a notion of variation of distributions using the Wasserstein distance, $W_r$, in the set $\mathcal{W}_r(\mathbb{R}^d)$ of probability measures with finite  {$r$th} moments,  {where} $r \geq 1 $. This notion generalizes the concept of  variance for random distributions over $\mathbb{R}^d$. This quantity can be defined as   
$$
V_r \left( \mu_1, \dots, \mu_J \right) = \inf_{\eta\in \mathcal{W}_r 
\left( \mathbb{R}^d\right)} \left\{\frac 1 J \sum_{j=1}^J {W}_r^r (\mu_j,\eta) \right\}^{1/r},
$$ 
which measures the spread of the distributions. Then, to measure closeness to a deformation model, we take a look at the minimal variation among warped distributions, a quantity that we could consider as a minimal alignment cost. Under some mild conditions, a deformation model holds if and only if this minimal alignment cost is null and we can base our assessment of a deformation model on this quantity. 

As in \cite{MR1704844, MR1625620}, we provide results (a Central Limit Theorem and bootstrap versions) that enable to reject that the minimal alignment cost exceeds some threshold, and hence to conclude that it is below that threshold. Our results are given in a setup of general, nonparametric classes of warping functions.  {We also provide results in the somewhat more restrictive setup where} one is interested in the more classical goodness-of-fit problem for the deformation model. Note that a general Central Limit Theorem is available for  {the} Wasserstein distance in \cite{2017arXiv170501299D}.

The paper is organized as follows. The main facts about Wasserstein variation are presented in Section~\ref{sec:2}, together with the key idea that fit to a deformation model can be recast in terms of the minimal Wasserstein variation among warped versions of the distributions. In Section~\ref{sec:3}, we prove some  {Lipschitz} bounds  for the law of empirical Wasserstein variations as well as of minimal alignment costs on $\mathbb{R}^d$.  {As a consequence of these results, the} quantiles of the minimal warped variation criterion can be consistently estimated by some suitable bootstrap quantiles, which can be approximated by simulation, yielding consistent tests of fit to deformation models, provided that the empirical criterion has  {a} regular limiting distribution. 

Central Limit Theorems for empirical minimal Wasserstein variation are further explored for univariate distributions in Section~\ref{sec:4},  {which covers} nonparametric deformation models, and in Section~\ref{sec:5},  {which presents} a sharper analysis for the case of semiparametric deformation models. These sections  {describe} consistent tests for deformation models in the corresponding setups. Section~\ref{sec:6}  {reports} some simulations  {assessing} the quality of the bootstrap procedure. Finally, proofs are  {gathered in the Appendix}.

\section{Wasserstein variation and deformation models for distributions\label{sec:2}}

Much recent work has been conducted to measure the spread or the inner structure of a collection of distributions. In this paper, we define a notion of variability which relies on  the notion of Fr\'echet mean for the space of {  probabilities} endowed with the Wasserstein metrics, of which we will recall the definition hereafter. First, for any integer $d\geq 1$,  consider  the set $\mathcal{W}_r ( \R^d )$ of probabilities with finite  {$r$th} moment. For $\mu$ and $\nu$ in $\mathcal{W}_r ( \R^d )$, we denote by $\Pi(\mu, \nu)$ the set of all probability measures $\pi$ over the product set $\R^d \times\R^d $ with first (respectively second) marginal $\mu$ (respectively $\nu$). The $L_r$ transportation cost between these two measures is defined as
\begin{equation*}
\label{eq:infwasser}
 {  W}_r(\mu, \nu)^r = \inf_{ \pi \in \Pi(\mu, \nu)} \int \left\Vert x - y\right\Vert ^r d \pi(x,y).
\end{equation*}
This transportation cost  {makes it possible to} endow the set $\mathcal{W}_r (\R^d )$ with the metric $W_r(\mu, \nu)$. More details on Wasserstein distances and their links with optimal transport problems can be found, e.g., in  \cite{rachev, villani2009optimal}.

Within this framework, we can define  a global measure of separation of a collection of probability measures as follows.  Given $\mu_1,\ldots,\mu_J\in\mathcal{W}_r(\mathbb{R}^d)$, let
$$
V_r\left( \mu_1, \dots, \mu_J \right) = \inf_{\eta\in \mathcal{W}_r(\mathbb{R}^d)} \left\{\frac 1 J \sum_{j=1}^J {W}_r^r (\mu_j,\eta) \right\}^{1/r}
$$
be the Wasserstein $r$-variation of $\mu_1,\ldots,\mu_J$ or the variance of the $\mu_j$s. 

The special case $r=2$ has been studied in the literature.  {The} existence of a minimizer of the map $\eta \mapsto  \{{W}_2^2 (\mu_1,\eta) + \cdots + W_2^2 (\mu_J, \eta)\}/J$ is proved in~\cite{agueh2010barycenters}, as well as  {its} uniqueness under some smoothness assumptions. Such a minimizer, $\mu_B$, is called a barycenter or Fr\'echet mean of $\mu_1,\ldots,\mu_J$. Hence,
$$
V_2\left( \mu_1, \dots, \mu_J \right) = \left\{\frac 1 J \sum_{j=1}^J {W}_2^2 (\mu_j,\mu_B) \right\}^{1/2}.
$$ 
Empirical versions of the barycenter are analyzed in~\cite{boissard2014, LGouic2016}. Similar ideas have also been developed in \cite{cuturi2014fast,Bigot12}.

This quantity, which is an extension of the variance for probability distributions is a good candidate 
to evaluate the concentration of a collection of measures around  {their} Fr\'echet mean. 
In particular, it can be used to measure  {the} fit to a  distribution deformation model.
More precisely, assume as in the Introduction that we observe $J$ independent  {random} samples
with sample $j \in \{ 1,\ldots,J\}$ consisting of iid observations $X_{1,j}, \ldots, X_{n,j}$ with common
distribution~$\mu_j$. We assume that $\mathcal{G}_j$ is a family (parametric or nonparametric)
of invertible warping functions and denote $\mathcal{G}=\mathcal{G}_1\times \cdots \times\mathcal{G}_J $.
Then, the deformation model assumes that
\begin{equation}\label{eq:H0}
\mbox{ there exists } (\varphi^ *_1, \dots , \varphi^*_J  )\in \mathcal{G} \mbox{ and iid } 
(\varepsilon_{i,j})_{1 \leq i \leq n, 1 \leq j\leq J} \mbox{ such that for all }  j \in \{1, \ldots, J\}, \, \\
X_{i,j} = (\varphi^*_j)^{-1} (\varepsilon_{i,j}).
\end{equation} 
Equivalently, the deformation model (\ref{eq:H0}) means that there  {exists} $(\varphi^ *_1, \dots , \varphi^*_J )\in \mathcal{G}$ such that the collection of $\varphi^*_j(X_{i,j})$s taken over all $j \in \{1, \ldots, J \}$ and $i \in \{ 1, \ldots, n\}$ is iid or, if we write $\mu_j(\varphi_j)$ for the distribution
of $\varphi_j(X_{i,j})$, that there exists $(\varphi^ *_1, \dots , \varphi^*_J) \in \mathcal{G}$ such that
$\mu_1(\varphi_1^*)=\cdots=\mu_J(\varphi_J^*)$.

We propose to use the Wasserstein variation to measure  {the} fit  {of} model~\eqref{eq:H0} through the minimal alignment cost
\begin{equation}\label{DefAr}
A_r(\mathcal{G}) = \inf_{(\varphi_1,\ldots,\varphi_J)\in \mathcal{G}}V_r^r\left\{ \mu_1(\varphi_1), \dots, \mu_J(\varphi_J) \right\}.
\end{equation}
Let us assume that $\mu_1(\varphi_1), \dots, \mu_J(\varphi_J)$, $(\varphi_1,\ldots,\varphi_J)\in \mathcal{G}$ are in $\mathcal{W}_r(\mathbb{R}^d)$. If the deformation model (\ref{eq:H0}) holds, then $A_r(\mathcal{G})=0$. Under the additional mild assumption that the minimum in (\ref{DefAr}) is attained, we have that the deformation model can be equivalently formulated as $A_r(\mathcal{G}) = 0$ and a goodness-of-fit test to the deformation model becomes, formally, a test of
\begin{equation}\label{classicalTest}
\mathcal{H}_0 :\, A_r(\mathcal{G})=0 \quad \mbox{vs.} \quad \mathcal{H}_a:\, A_r(\mathcal{G})>0.
\end{equation}
A testing procedure can be based on the empirical version of $A_r(\mathcal{G})$, namely,
\begin{equation}\label{DefArEmp}
A_{n,r}(\mathcal{G}) = \inf_{(\varphi_1,\ldots,\varphi_J)\in \mathcal{G}}V_r^r\left\{ \mu_{n,1}(\varphi_1), \dots, \mu_{n,J}(\varphi_J) \right\},
\end{equation}
where $\mu_{n,j}(\varphi_j)$ denotes the empirical measure on $\varphi_j(X_{1,j}),\ldots,\varphi_j(X_{n,j})$. We would reject the deformation model (\ref{eq:H0}) for large values of $A_{n,r}(\mathcal{G})$.

As noted in \cite{MR1704844,MR1625620}, the testing problem (\ref{classicalTest}) can be considered as a mere  sanity check for the deformation model, since lack of rejection of the null does not provide statistical evidence that the deformation model holds. Consequently, as in the cited references, we will also consider the alternative testing  problem
\begin{equation}\label{eq:H0_test_nonpar} 
\mathcal{H}_0 :\, A_r(\mathcal{G})\geq \Delta_0 \quad \mbox{vs.} \quad \mathcal{H}_a:\, A_r(\mathcal{G})<\Delta_0,
\end{equation}
where $\Delta_0 >0$ is a fixed threshold. With this formulation the  test decision of rejecting  the null hypothesis implies that there is statistical evidence that the deformation model is approximately true. In this case, rejection would correspond to small observed values of $A_{n,r}(\mathcal{G})$. In  {subsequent} sections, we provide theoretical results that allow the computation of approximate critical values and $p$-values for the testing problems (\ref{classicalTest}) and (\ref{eq:H0_test_nonpar}) under suitable assumptions. 

\section{Bootstraping Wasserstein's variations\label{sec:3}}

We present now some general results on Wasserstein distances that will be  applied to estimate the asymptotic distribution of the minimal alignment cost statistic, $A_{n,r}(\mathcal{G})$, defined in (\ref{DefArEmp}). In this section, we write $\mathcal{L}(Z)$ for the law of any random variable $Z$. We note the abuse of notation in the following, in which ${W}_r$ is used both for the Wasserstein distance on $\mathbb{R}$ and on $\mathbb{R}^d$, but this should not cause much confusion.

Our first result shows that the laws of empirical transportation costs are
continuous (and even Lipschitz) functions of the underlying distributions.
\begin{theo}\label{prop:transportationcost}
Set $\nu, \nu', \eta$  probability measures in $\mathcal{W}_r (\mathbb{R}^d)$, $Y_1,\ldots,Y_n$  
iid random vectors with common law $\nu$, $Y'_1,\ldots,Y'_n$, iid with law  $\nu'$ and
write $\nu_n$, $\nu_n'$ for the corresponding empirical measures.  Then
$${W}_r[\mathcal{L}\{{W}_r(\nu_n,\eta)\},\mathcal{L}\{{W}_r(\nu'_n,\eta)\}]\leq {W}_r(\nu,\nu').$$
\end{theo}
The deformation assessment criterion introduced in Section~\ref{sec:2} is based on the  
{Wasserstein $r$-variation} of distributions, $V_r$. It is convenient to note that $V_r^r(\nu_1,\ldots,\nu_J)$
can also be expressed as
\begin{equation}
\label{eq:new_formulation}
V_r^r(\nu_1,\ldots,\nu_J)=\inf_{\pi\in \Pi(\nu_1,\ldots,\nu_J)} \int T(y_1,\ldots,y_J)d\pi(y_1,\ldots,y_J),
\end{equation}
where $\Pi(\nu_1,\ldots,\nu_J)$ denotes the set of probability measures on $\mathbb{R}^d$ with marginals
$\nu_1,\ldots,\nu_J$ and 
$$ 
T(y_1,\ldots,y_J)=\min_{z\in\mathbb{R}^d} \frac 1 J \sum_{j=1}^J \|y_j-z\|^r.
$$ 

Here we are interested in empirical Wasserstein $r$-variations, namely, the $r$-variations computed from the empirical measures $\nu_{n_j,j}$ coming from independent samples $Y_{1, j},\ldots,Y_{n_j, j}$ of  {iid} random variables with distribution $\nu_j$. Note that in this case, problem \eqref{eq:new_formulation} is a linear optimization problem  for which a minimizer always exists.

As before, we consider the continuity of the law of empirical Wasserstein $r$-variations with respect to the underlying probabilities. This is covered in the next result.

\begin{theo}\label{prop:rvariation}
With the above notation,
$$
{W}_r^r
[ \mathcal{L} \{ V_r(\nu_{n_1,1},\ldots,\nu_{n_J,J})\}, \mathcal{L} \{ V_r(\nu'_{n_1,1},\ldots,\nu'_{n_J,J})\}]
\leq \frac 1 J \sum_{j=1}^J {W}_r^r(\nu_j,\nu'_j).
$$
\end{theo}
A useful consequence of the above results is that empirical Wasserstein distances or $r$-variations can be bootstrapped under rather general conditions. To be more precise, in Theorem~\ref{prop:transportationcost} we take $\nu'=\nu_n$, the empirical measure on $Y_1,\ldots,Y_n$, and consider a bootstrap sample $Y_1^*,\ldots,Y_{m_n}^*$ of iid (conditionally given $Y_1,\ldots,Y_n$) observations with common law $\nu_n$. We will assume that  the resampling size $m_n$ satisfies $m_n\to \infty$, $m_n=o(n)$ and write $\nu_{m_n}^*$ for the empirical measure on $Y_1^*,\ldots,Y_{m_n}^*$ and $\mathcal{L}^*(Z)$ for the conditional law of $Z$ given $Y_1,\ldots,Y_n$. Theorem~\ref{prop:transportationcost} now reads 
$$
{W}_r [ \mathcal{L}^* \{ {W}_r(\nu^*_{m_n},\nu)\}, \mathcal{L}({W}_r \{ \nu_{m_n},\nu) \} ]\leq {W}_r(\nu_n,\nu).
$$
Hence, if ${W}_r(\nu_n,\nu)=O_{\p}(1/r_n)$ for some sequence $r_n>0$ such that $r_{m_n}/r_n\to 0$ as $n\to\infty$, then, using the fact that ${W}_r \{ \mathcal{L}(aX),\mathcal{L}(aY) \} = a{W}_r \{ \mathcal{L}(X),\mathcal{L}(Y)\}$ for $a>0$, we see that 
\begin{equation*}
{W}_r[\mathcal{L}^*\{r_{m_n}{W}_r(\nu^*_{m_n},\nu)\},\mathcal{L}\{r_{m_n}{W}_r(\nu_{m_n},\nu)\}]\leq 
\frac{r_{m_n}}{r_n} \, r_n {W}_r(\nu_n,\nu)\to 0
\end{equation*}
 in probability.

 {Assume} that, in addition,   $r_n{W}_r(\nu_n,\nu)\rightsquigarrow \gamma \left(\nu \right)$ 
for a smooth distribution $\gamma \left(\nu \right)$. If $\hat{c}_n(\alpha)$ denotes the $\alpha$th quantile of the conditional distribution $\mathcal{L}^*\{r_{m_n}{W}_r(\nu^*_{m_n},\nu)\}$, then
\begin{equation}\label{eq:bootquantile}
\lim_{n\to \infty} \p\{ r_n{W}_r(\nu_n,\nu) \leq \hat{c}_n(\alpha)\} = \alpha;
\end{equation}
see, e.g., Lemma 1 in \cite{JanssenPauls}. 
We conclude in this case that  the quantiles of $r_n{W}_r(\nu_n,\nu)$
can be consistently estimated by the bootstrap quantiles, $\hat{c}_n(\alpha)$, which, in turn, can be approximated through 
Monte Carlo simulation.

As an example, if $d=1$ and $r=2$, under integrability and smoothness assumptions on $\nu$, we have 
$$
\sqrt{n}\, {W}_2(\nu_n,\nu)\rightsquigarrow \left[ \int_{0}^1 \frac{B^2(t)}{f^2\{ F^{-1}(t)\}} dt\right]^{1/2},
$$ 
as $n \to \infty$, where $f$ and $F^{-1}$ are the density and the quantile function of $\nu$,  {respectively}; see~\cite{delBarrioGineUtzet}).  {Therefore, Eq.}~(\ref{eq:bootquantile}) holds. Bootstrap results have also been provided in \cite{ebert}.

For the  deformation model (\ref{eq:H0}),   statistical inference  is based on 
$A_{n,r}(\mathcal{G})$, introduced in (\ref{DefArEmp}).
Now consider $A'_{n,r}(\mathcal{G})$, the corresponding version obtained from samples with underlying distributions $\mu_j'$. Then, a version of Theorem~\ref{prop:rvariation} is valid for these minimal alignment costs, provided  {that} the deformation classes are uniformly Lipschitz, namely, under the assumption that, for all $j \in \{ 1, \ldots, J\}$, 
\begin{equation}\label{UnifLipsch}
L_j = \sup_{x\ne y, \varphi_j\in \mathcal{G}_j} \frac{\|\varphi_j(x)-\varphi_j(x)\|}{\|x-y\|}
\end{equation}
 {is} finite.

\begin{theo} \label{th:bootstrap}
If $L=\max(L_1,\ldots,L_j)<\infty$, with $L_j$ as in~\eqref{UnifLipsch}, then  
$${W}_r^r\left[\mathcal{L}[\{A_{n,r}(\mathcal{G})\}^{1/r}],\mathcal{L}[\{A'_{n,r}(\mathcal{G})\}^{1/r} ] \right]
\leq L^r \frac 1 J \sum_{j=1}^J {W}_r^r (\mu_j,\mu'_j).$$
\end{theo}
Hence, the Wasserstein distance of the variance of two collections of distributions can be controlled using the distance between the distributions. The main consequence of this fact is that the minimal alignment cost can  {also be} bootstrapped as soon as a distributional limit theorem exists for $A_{n,r}(\mathcal{G})$, as in the discussion above. In Sections~\ref{sec:4} and \ref{sec:5} below, we present distributional results of this type in the  {one-dimensional} case. We note that, while general Central Limit Theorems for the empirical transportation cost are not available in dimension $d>1$, some recent progress has been made  {in this direction}; see, e.g.,~\cite{rippl15} for Gaussian distributions and~\cite{SommerfeldMunk}, which gives  such  results for distributions on $\mathbb{R}^d$ with finite support. Further advances  {along these lines} would  {make it possible} to extend the results in the following section to higher dimensions.

\section{Assessing  {the} fit  {of} nonparametric deformation models\label{sec:4}}

 {In this section and subsequent ones, we focus} on the case $d=1$ and $r=2$.  {Thus we} will simply write $A(\mathcal{G})$ and $A_n(\mathcal{G})$ (instead of $A_2(\mathcal{G})$ and $A_{2,n}(\mathcal{G})$) for the minimal alignment cost and its empirical version, defined in (\ref{DefAr}) and (\ref{DefArEmp}). Otherwise we keep the notation  {from} Section~\ref{sec:2}, with $X_{1,j},\ldots,X_{n,j}$ iid random variables with law $\mu_j$ being one of the $J$ independent samples. Now $\mathcal{G}_j$ is a class of invertible warping functions from $\mathbb{R}$ to $\mathbb{R}$, which we assume to be increasing. We note that in this case the barycenter of a set of  {probability measures} $\mu_1,\ldots,\mu_J$ with distribution functions $F_1,\ldots,F_J$ is the probability having quantile function $F_B^{-1} = (F_1^{-1} + \cdots + F_J^{-1})/J$; see, e.g.,~\cite{agueh2010barycenters}. We observe further that $\mu_j(\varphi_j)$ is determined by the quantile
function $\varphi_j\circ F_j^{-1}$. We will write 
\begin{equation*}\label{barycenternotation}
F_{B}^{-1}(\varphi)=\frac 1 J \sum_{j=1}^J \varphi_j\circ F_j^{-1}
\end{equation*}
for the quantile function of the barycenter of $\mu_1(\varphi_1),\ldots,\mu_J(\varphi_J)$, while $\rightsquigarrow$
will denote convergence in distribution.

In order to prove a Central Limit Theorem for $A_n(\mathcal{G})$, we need to make assumptions on the integrability
and regularity of the distributions $\mu_1, \ldots, \mu_J$ as well as on the smoothness of the warping functions.
We consider first the assumptions on the distributions. For each $j \in \{ 1,\ldots,J\}$,  {the distribution function associated with $\mu_j$ is denoted} $F_j$. We will assume that $\mu_j$ is supported 
on a (possibly unbounded) interval in the interior of which $F_j$ is $C^2$ and $F'_j=f_j>0$ and satisfies
\begin{equation}
\label{hyp:tassio13} 
\sup_{x}  \frac{F_j(x) \{ 1- F_j(x) \} { |f'_j(x)|}}{f_j(x)^ 2} < \infty,
\end{equation}
and, further, that for some $q>1$,
\begin{equation}
\label{hyp:Rajput}
\displaystyle
\int_0 ^1 \frac{\left\{t \left(1-t\right)\right\}^{q/2}}{ [f_j \{ F_j^{-1} ( t  )  \} ] ^q}dt <\infty
\end{equation}
and for some $r>4$,
\begin{equation}
\label{hyp:tassio12}
\E  ( \vert X_j \vert^ r ) < \infty.
\end{equation}

Assumption (\ref{hyp:tassio13}) is a classical regularity requirement for the use of strong approximations for the quantile process, as in
\cite{csorgohorvath93,delBarrioGineUtzet}. Our proof relies on the use of these techniques. As for  
(\ref{hyp:Rajput}) and (\ref{hyp:tassio12}), they are mild integrability conditions. If $F_j$ has regularly varying tails of order $-r$ (e.g.,
Pareto tails) then both conditions hold --- and also (\ref{hyp:tassio13}) --- so long as $r>4$ and $1<q<2r/(r+2)$. 
Of course the conditions are fulfilled by distributions with lighter tails such as exponential
or Gaussian laws for any $q\in(1,2)$. 

Turning to the assumptions on the classes of warping functions, we recall that a uniform  {Lipschitz}
condition was needed for the approximation bound in Theorem~\ref{th:bootstrap}. For the Central Limit 
Theorem in this section, we
need some refinement of that condition, the extent of which will depend on the integrability exponent $q$ in 
(\ref{hyp:Rajput}), as follows. We set \mbox{$p_0 = \max \{{q}/{(q-1)}, 2 \}$} and define on 
$\mathcal{H}_j=C^1 (\mathbb{R}) \cap  L^{p_0}( X_j ) $ the norm
$$
\Vert h_j  \Vert_{\mathcal{H}_j}  = \sup_{} \vert h_j'(x) \vert + 
\E  \{  \vert h_j(X_j)  \vert ^{p_0} \}^{1/p_0},
$$ 
and on the product space $\mathcal{H}_1 \times \cdots \times\mathcal{H}_J$, $\Vert h \Vert_{\mathcal{H}}  
=  \Vert h_1 \Vert_{\mathcal{H}_1} + \cdots + \Vert h_J \Vert_{\mathcal{H}_J}$. We further assume that
\begin{equation}
\label{hyp:regul}
\mathcal{G}_j\subset \mathcal{H}_j \mbox{ is compact for } \Vert \cdot \Vert_{\mathcal{H}_j}
\mbox{ and } \sup_{h  \in \mathcal{G}_j}  \vert h^\prime (x^h _n)- h^\prime(x) 
 \vert \underset{ \sup_{h  \in \mathcal{G}_j} \vert x_n^h- x  \vert \rightarrow 0}\longrightarrow 0,
\end{equation}
and, finally, that  for some $r>\mbox{max}(4,p_0)$,
\begin{equation}
\label{hyp:tassio12theta}
\E \sup_{h \in \mathcal{G}_j}  \vert h ( X_j) \vert ^ r <\infty.
\end{equation}

We note that (\ref{hyp:tassio12theta}) is a slight strengthening of the uniform moment bound already contained
in (\ref{hyp:regul}); we could ta\-ke~$p_0>\max\{{q}/{(q-1)},4\}$ in~(\ref{hyp:regul}) and~(\ref{hyp:tassio12theta}) would follow. 

Our next result gives a Central Limit Theorem for $A_n(\mathcal{G})$ under the assumptions on the distributions and deformation classes described above. The limit can be simply described in terms of a centered Gaussian process indexed by the set of minimizers of the
variation functional, namely,
$U(\varphi)=V_2^2\{ \mu_1(\varphi_1), \dots, \mu_J(\varphi_J)\}.$
An elementary computation shows that 
$$
\{U^{1/2}(\varphi)-U^{1/2}(\tilde{\varphi}) \}^2\leq \frac 1 J\sum_{j=1}^J \E \{\varphi_j(X_j)-\tilde{\varphi}_j(X_j) \}^2,
$$
from which we conclude continuity of $U$ with respect to $\|\cdot\|_{\mathcal{H}}$. In particular, the set
\begin{equation*}\label{gammaset}
\Gamma=\Big\{\varphi\in\mathcal{G}:\, U(\varphi)=\inf_{\phi \in \mathcal{G}}U(\phi)\Big\}
\end{equation*}
is a nonempty compact subset of $\mathcal{G}$.

\begin{theo}
\label{th:testgen}
Assume that $B_1, \ldots, B_J$ are mutually independent Brownian bridges. Set
$$
c_{j}(\varphi) =2 \int_0^1  \varphi_j^\prime \circ F_j^{-1}  \{\varphi_j\circ F_j^{-1}-F_B^{-1}(\varphi) \} \frac{B_{j}}{f_j\circ F_j^{-1}}
$$
and $C(\varphi) = \{c_1(\varphi) + \cdots + c_J(\varphi)\}/J$, $\varphi\in\mathcal{G}$. Then, under  assumptions (\ref{hyp:tassio13})--(\ref{hyp:tassio12theta}), $C$ is a centered Gaussian process on $\mathcal{G}$ with trajectories that are almost surely continuous with respect to $\|\cdot\|_\mathcal{H}$. Furthermore,
\begin{eqnarray*}
\sqrt{n}\left\{A_n(\mathcal{G})-A(\mathcal{G})\right\} \rightsquigarrow \min_{\varphi  \in \Gamma} C(\varphi).
\end{eqnarray*}
\end{theo}

A proof of Theorem~\ref{th:testgen} is given in the Appendix. The random variables 
$$
\int_0^1  \varphi_j^\prime\circ F_j^{-1}\frac{B_j}{ {f_j\circ F_j^{-1}}} \{\varphi_j \circ F_j^{-1} -F^{-1}_B(\varphi) \}
$$
are centered Gaussian, with variance 
\begin{equation*}
\int_{[0,1]^2} 
\{\min(s,t)- st\}  
\frac{\varphi_{j}^{\prime}\{F_j^{-1}(t)\}}{{f_j\{F_j^{-1}(t)\}}}
\left[\varphi_{j}\{ F_j^{-1}(t)\} -F^{-1}_B(\varphi)(t)\right]
\frac{\varphi_{j}^{\prime}\{F_j^{-1}(s)\}}{{f_j\{F_j^{-1}(s)\}}} 
\left[\varphi_{j}\{ F_j^{-1}(s)\} - F^{-1}_B(\varphi)(s)\right] ds dt.
\end{equation*}
In particular, if $U$ has a unique minimizer the limiting distribution in Theorem~\ref{th:testgen} is  {Gaussian}. However, our result works in more generality, even without uniqueness assumptions.

We remark also that although we have focused for simplicity on the case of samples of equal size, the case of different sample sizes $n_1,\ldots,n_J$ can also be handled with straightforward changes. More precisely, let us  write $A_{n_1,\ldots,n_J}(\mathcal{G})$ for the minimal alignment cost computed from the empirical distribution of the samples and assume that $n_j \rightarrow \infty$ and
$$
\frac{n_j}{n_1 + \cdots + n_J} \rightarrow ( \gamma_j ) ^2>0,
$$ 
then with straightforward changes in our proof we can see that 
\begin{equation*}\label{eq:condi_ndif}
\sqrt{\frac{n_1 \dots n_J}{\left(n_1 + \dots + n_J\right)^{J-1}}}\,  
\left\{ A_{n_1,\ldots,n_J}(\mathcal{G})-A(\mathcal{G})\right\} \rightsquigarrow  \min_{\varphi\in\Gamma} \tilde{C}(\varphi),
\end{equation*} 
where $\tilde{C}(\varphi)= \{ \tilde{c}_1(\varphi) + \cdots + \tilde{c}_J(\varphi)\}/J$ and $\tilde{c}_j(\varphi)=\big(\Pi_{p\neq j} \gamma_p\big)c_j(\varphi)$.

If we try, as argued in Section~\ref{sec:2}, to base  {on $A_n(\mathcal{G})$} our assessment of fit to the deformation model (\ref{eq:H0}), we should note that the limiting distribution in Theorem~\ref{th:testgen} depends on the unknown distributions $\mu_j$ and cannot be used for the computation of approximate critical values or $p$-values without further adjustments. We show now how this can be done in the case of the testing problem (\ref{eq:H0_test_nonpar}), namely, the test of
$$
\mathcal{H}_0 :\, A_r(\mathcal{G})\geq \Delta_0 \quad \mbox{vs.} \quad \mathcal{H}_a:\, A_r(\mathcal{G})<\Delta_0,
$$
for some fixed threshold $\Delta_0>0$, through the use of a bootstrap procedure.

Let us consider bootstrap samples $X^*_{1,j},\ldots,X^*_{m_{n},j}$ 
of iid observations sampled from $\mu_{n,j}$, the empirical distribution on $X_{1,j},\ldots,X_{n,j}$.
We write $\mu^*_{m_n,j}$ for the empirical measure on
$X^*_{1,j},\ldots,X^*_{m_n,j}$ and introduce
$$
A^{*}_{m_n} (\mathcal{G})=\inf_{\varphi\in\mathcal{G}} V^2_2 \{\mu^*_{m_n,1}(\varphi_1 ),\ldots,\mu^*_{m_n,J}(\varphi_J) \}.
$$ 
Now, we base our testing procedure on the  conditional $\alpha$-quantiles (given the $X_{i,j}$s) of $\sqrt{m_n} \, \{A^{*}_{m_n} (\mathcal{G})-\Delta_0\}$,
which we denote $\hat{c}_n(\alpha;\Delta_0)$. Our next result, which follows from Theorems \ref{th:bootstrap} and \ref{th:testgen},  
shows that the test that rejects $\mathcal{H}_0 $ when
$$
\sqrt{n} \, \{A_n(\mathcal{G})-\Delta_0\} < \hat{c}_n(\alpha;\Delta_0)
$$
is a consistent test of approximate level $\alpha$ for (\ref{eq:H0_test_nonpar}).
We note that the bootstrap quantiles $\hat{c}_n(\alpha; \Delta_0)$ can be computed using Monte Carlo simulation.
\begin{corollary}\label{cor:bootstrap}
If $m_n\to \infty$, and $m_n=O(\sqrt{n})$, then under assumptions~(\ref{hyp:tassio13})--(\ref{hyp:tassio12theta})
\begin{equation*}
\label{eq:res_cor}
\p\big[\sqrt{n}\{A_n(\mathcal{G})-\Delta_0\} < \hat{c}_n(\alpha;\Delta_0) \big] \to 
\begin{cases}
0 &\mbox{if } A({\mathcal G})>\Delta_0,\\
\alpha & \mbox{if } A({\mathcal G})=\Delta_0,\\
1 & \mbox{if } A({\mathcal G}) < \Delta_0.
\end{cases}
\end{equation*}
\end{corollary}

Rejection in the testing problem \eqref{eq:H0_test_nonpar} would result, as noted in Section~\ref{sec:2}, in statistical evidence supporting that the deformation model holds approximately, and hence that related registration methods can be safely applied. If, nevertheless, we were interested in gathering statistical evidence against the deformation model, then we should consider the classical goodness-of-fit problem (\ref{classicalTest}). Some technical difficulties arise then. Note that if the deformation model holds, that is, if $A(\mathcal{G})=0$, then we have $\varphi_j\circ F_j^{-1}=F_B^{-1}(\varphi)$ for each $\varphi\in\Gamma$, which implies that the result of Theorem~\ref{th:testgen} becomes $\sqrt{n}\, A_n(\mathcal{G}) \rightsquigarrow 0$. Hence, a nondegenerate limit law for $ A_n(\mathcal{G})$ in this case requires a more refined analysis, that we handle in the next section.

\section{Goodness-of-fit in semiparametric deformation models\label{sec:5}}

In many cases, deformation functions can be made more specific in the sense that they follow a known shape depending on parameters that may differ  for sample to sample. In our approach to the classical goodness-of-fit problem (\ref{classicalTest}), we consider a parametric model in which $\varphi_j=\varphi_{\theta_j}$ for some  {finite-dimensional} parameter $\theta_j$ that describes the warping effect within a fixed shape. Now, that the deformation model holds means that there exist $\theta^*=(\theta_1^*,\dots,\theta_J^*)$ such that for all $i \in \{ 1, \ldots, n\}$ and $j \in \{ 1, \ldots, J\}$, $X_{i,j}= \varphi_ {\theta^*_j}^{-1} (\varepsilon_{i,j} )$. Hence, from now on, we will consider the following family of deformations, indexed by a parameter $ \lambda \in \Lambda \subset \R^p$: 
$$ 
\varphi :  \Lambda \times \mathbb{R}  \rightarrow  \mathbb{R} : 
\left(\lambda,x\right) \mapsto  \varphi_{\lambda} (x).
$$

The classes $\mathcal{G}_j$ become now $\{\varphi_{\theta_j}:\, \theta_j\in\Lambda\}$. We denote $\Theta= \Lambda^J$ and write $A_n(\Theta)$ and $A(\Theta)$ instead of $A_n(\mathcal{G})$ and $A(\mathcal{G})$. We also use the simplified notation $\mu_j(\theta_j)$ instead of $\mu_j ( \varphi_{\theta_j} )$, $F_B \left( {\theta} \right)$ for $F_B \left( \varphi_{\theta_1}, \dots, \varphi_{\theta_J} \right)$ and similarly for the empirical versions. Our main goal is to  {identify} a weak limit theorem for $A_n(\Theta)$ under the null in (\ref{classicalTest}). Therefore, throughout this section, we assume that model~(\ref{eq:H0}) holds. This means, in particular, that the quantile functions of the samples satisfy $F_j^{-1}=\varphi_{\theta_j^*}^{-1}\circ G^{-1}$, with $G$ the distribution function of the $\varepsilon_{i,j}$s. As before, we assume that the warping functions are invertible and increasing, which now means that, for each $\lambda \in \Lambda$, $\varphi_\lambda$  is an invertible, increasing function. 

It is convenient at this point to introduce the notation
\begin{equation*}\label{psifucntion}
\psi_j(\lambda,x)=\varphi_\lambda\{\varphi_{\theta_j^*}^{-1}(x)\}
\end{equation*}
for all $j \in \{ 1, \ldots, J\}$ and $\varepsilon$ for a random variable with the same distribution as the $\varepsilon_{i,j}$. Note that $\psi_j(\theta_j^*,x)=x$.

Now, under smoothness assumptions on the functions $\psi_j$ that we present in detail below, if the parameter space is compact then the function
$$
U_n(\theta_1,\ldots,\theta_J) = V_2^2 \{ \mu_{n,1}(\theta_1),\ldots,\mu_{n,J}(\theta_J)\}
$$
admits a minimizer that we will denote by $\hat{\theta}_n$, i.e.,
\begin{equation*}\label{argminemp}
\hat{\theta}_n\in \underset {\theta \in\Theta} {\mbox{argmin}} \, U_n(\theta).
\end{equation*}
Of course, since we are assuming that the deformation model holds, we know that $\theta^*$ is a minimizer of 
$$U(\theta_1,\ldots,\theta_J)=V_2^2 \{\mu_{1}(\theta_1),\ldots,\mu_{J}(\theta_J)\}.$$ 
For a closer analysis of the asymptotic behavior of $A_n(\Theta)$ under the deformation model, we need to make the following identifiability assumption
\begin{equation}
\label{hyp:identifTassio}
\theta^* \text{ belongs to the interior of }\Lambda \mbox{ and is the unique minimizer of } U.
\end{equation}
Note that, equivalently, this means that $\theta^*$ is the unique zero of $U$.

As in the case of nonparametric deformation models, we need to impose some conditions on the class of warping functions and on the distribution of the errors, the $\varepsilon_{i,j}$s. For the former, we write $D$ or $D_u$ for derivative operators with respect to parameters. Hence, for instance, $D\psi_j(\lambda,x)=(D_1\psi_j(\lambda,x), \ldots,D_p\psi_j(\lambda,x))^{\top}$ is the vector consisting of partial derivatives of $\psi_j$ with respect to its first $p$ arguments evaluated at $(\lambda, x)$; $D^2\psi_j(\lambda,x)=(D_{u,v}\psi_j(\lambda,x))_{u,v}$ is the hessian matrix for fixed $x$ and so on.  {In what follows,} $\psi_j'(\lambda,x)$ and similar notation will stand for derivatives with respect to $x$. Then we will assume that  for each $j \in \{ 1,\ldots,J\}$, $u,v \in \{ 1,\ldots,p\}$, and some $r>4$
\begin{equation}\label{jointsmoothness}
\psi_j(\cdot,\cdot) \mbox{ is } C^2,
\end{equation}
\begin{equation}\label{uniformmoments}
\E \Big\{ \sup_{\lambda \in \Lambda } \vert \psi_j(\lambda,\varepsilon) \vert^r \Big\} < \infty, \quad 
\E \Big\{ \sup_{\lambda \in \Lambda } \vert D_{u}\psi_j(\lambda,\varepsilon) \vert^ r \Big\} < \infty,\quad 
\E \Big\{ \sup_{\lambda \in \Lambda } \vert D_{u,v}\psi_j(\lambda,\varepsilon) \vert^ r \Big\} < \infty, 
\end{equation}
and
\begin{equation}\label{hyp:regul_th}
\psi_j'(\cdot,\cdot) \mbox{ is bounded on } \Lambda \times \mathbb{R} \; \mbox{ and } \;
\sup_{\lambda \in \Lambda}  \left\vert \psi_j'(\lambda,x^\lambda _n)-  
\psi_j'(\lambda,x)\right\vert \xrightarrow{ \sup_{\lambda \in \Lambda} \vert x_n^\lambda- x  \vert \rightarrow 0}0.
\end{equation}

Turning to the distribution of the errors, we will assume that
$G$ is $C^2$ with  $G'(x)=g(x)>0$  on some interval and
\begin{equation}
\label{hyp:tassio13_eps}  
\sup_{x} \frac{G(x)\left\{ 1- G(x)\right\} {  |g'(x)|}}{g(x)^ 2} < \infty.
\end{equation}
Additionally (but see the comments after Theorem~\ref{prop:tcl} below), we make the assumption that
\begin{equation}\label{hyp:tassioA3} 
\int_0 ^ 1 \frac{t(1-t)}{g^2\left\{G^ {-1} (t) \right\}}dt <\infty.
\end{equation}
Finally, before stating the asymptotic result for $A_n(\Theta)$, we introduce the $p\times p$ matrices
$$\Sigma_{i,i}=\frac{2(J-1)}{J^2}\int_0^1 D_i \psi_i \{\theta_i^*, G^{-1}(t)\} \psi_i \{\theta_i^*, G^{-1}(t)\}^{\top} dt,$$
$$\Sigma_{i,j}=-\frac{2}{J^2}\int_0^1 D_i \psi_i \{\theta_i^*, G^{-1}(t)\} \psi_i \{\theta_j^*, G^{-1}(t)\}^{\top} dt,\quad i\ne j$$
and the $pJ \times pJ$ matrix
\begin{equation*}\label{def:phi}
\Sigma=\left[ 
\begin{matrix}
\Sigma_{1,1} & \cdots & \Sigma_{1,J} \\
\vdots & & \vdots\\
\Sigma_{J,1} & \cdots & \Sigma_{J,J} 
\end{matrix}
\right].
\end{equation*}
 {This matrix} is symmetric and positive semidefinite. To see this, consider $x_1,\ldots,x_J \in\mathbb{R}^p$ and $x^{\top}=[x^{\top}_1,\ldots,x_J^{\top}]$. Note that
\begin{align*}
x'\Sigma x &= \frac 2{J^2}\int_0^1 
\bigg[\sum_i (J-1) \left[x_i \times D_i \psi_i \big\{\theta_i^*, G^{-1}(t) \big\}\right]^2 
- 2\sum_{i<j} \left[ x_i \times D_i \psi_i \big\{\theta_i^*, G^{-1}(t) \big\} \right] \left[ x_j \times D_j \psi_j \big\{\theta_j^*, G^{-1}(t) \big\} \right] \bigg] dt\\
&= \frac 2{J^2}\int_0^1 
\sum_{i<j} \left[ \big[x_i \times D_i \psi_i \big\{ \theta_i^*, G^{-1}(t) \big\} \big] - \big[ x_j \times D_j \psi_j \big\{ \theta_j^*, G^{-1}(t) \big\} \big] \right]^2 dt\geq 0.
\end{align*}
In fact, $\Sigma$ is positive definite, hence invertible, apart from some degenerate cases, For instance, if $p=1$,
$\Sigma$ is invertible unless all the functions $D_i \psi_i\{\theta_i^*, G^{-1}(t)\}$ are proportional.

We are ready now for the announced distributional limit theorem.
\begin{theo}
\label{prop:tcl}
Assume that the deformation model holds. Under  assumptions \eqref{hyp:identifTassio}--\eqref{hyp:tassio13_eps} 
$\hat{\theta}_n\to \theta^*$ 
in probability. If, in addition, $\Phi$ is invertible, then
$\sqrt{n} \, (\hat{\theta}_n - \theta^*) \rightsquigarrow \Sigma^{-1} Y$,
where $Y=(Y_1^{\top},\ldots,Y_J^{\top})^{\top}$ with 
$$Y_{j}=\frac{2}{J} \int_0^1 D \psi_j \{\theta_j^*,G^{-1}(t)\} \frac{\tilde{B}_j(t)}{g \{G^{-1}(t) \} }dt,$$
$\tilde{B_j}=B_j- (B_1 + \cdots + B_J)/J$ and $B_1, \ldots, B_J$ mutually independent Brownian bridges. Furthermore, if
\eqref{hyp:tassioA3} also holds, then
\begin{align*}
n A_n(\Theta) \rightsquigarrow  
\frac 1 J\sum_{j=1}^J \int_0^1 \bigg(\frac{\tilde{B}_j}{g \circ G^{-1}  }\bigg)^2 -\frac 1 2 Y^{\top}\Sigma^{-1}Y.
\end{align*}
\end{theo}   

We make a number of comments here. First, we note that, while, for simplicity, we have formulated Theorem~\ref{prop:tcl} assuming that the deformation model holds, the Central Limit Theorem for $\hat{\theta}_n$ still holds (with some additional assumptions and changes in $\Phi$) in the case when the model is false and $\theta^*$ is not the {true} parameter, but the one that gives the best (but imperfect) alignment.  {Given that} our focus here is the assessment of the deformation models, we refrain from pursuing this issue.

Our second comment is about the  {identifiability} condition (\ref{hyp:identifTassio}). At first sight it can seem to be too strong to be realistic. Actually, for some deformation models it could happen that $\varphi_\theta \circ \varphi_\eta=\varphi_{\theta*\eta}$ for some $\theta*\eta\in\Theta$. In this case, if $X_{i,j}=\varphi_{\theta_j^*}^{-1}(\varepsilon_{i,j})$ with $\varepsilon_{i,j}$ iid, then, for any $\theta$, $X_{i,j}=\varphi_{\theta*\theta_j^*}^{-1}(\tilde{\varepsilon}_{i,j})$ with $\tilde{\varepsilon}_{i,j}= \varphi_{\theta}({\varepsilon}_{i,j})$ which are also iid;  consequently, $(\theta* \theta^*_1,\ldots,\theta* \theta^*_J)$ is also a zero of $U$. This applies, for instance, to location and scale models. 

A simple fix to this issue is to select one of the signals as the reference, say the  {$J$th} signal, and assume that $\theta_J^*$ is known since it can be, in fact, chosen arbitrarily. The criterion function becomes then $\tilde{U}(\theta_1,\ldots,\theta_{J-1})=U(\theta_1,\ldots,\theta_{J-1},\theta_J^*)$. One could then make the (more realistic) assumption that $\tilde{\theta}^*=(\theta_1^*,\ldots,\theta_{J-1}^*)$ is the unique zero of $\tilde{U}$ and base the analysis on $\tilde{U}_n(\theta_1,\ldots,\theta_{J-1})=U_n(\theta_1,\ldots,\theta_{J-1},\theta_J^*)$ and $\hat{\tilde{\theta}}_n=\arg\min_{\tilde{\theta}} \tilde{U}_n(\tilde{\theta})$. The results in this section can be adapted almost verbatim to this setup. In particular, 
$$
\sqrt{n} \, (\hat{\tilde{\theta}}_n-{\tilde{\theta}^*})
\rightsquigarrow \tilde{\Sigma}^{-1} \tilde{Y}, 
$$
with $\tilde{Y}^{\top}=(Y_1^{\top},\ldots,Y^{\top}_{J-1})$ and  $\tilde{\Sigma} = (\Sigma_{i,j})_{1\leq i,j\leq J-1}$. Again, the invertibility of $\tilde{\Sigma}$ is almost granted. In fact, arguing as above, we see that and $\tilde{\Sigma}$ is positive definite if the function $D \psi_i\{\theta_i^*,G^{-1}(t)\}$ is not null for all $i \in \{ 1, \ldots, J-1\}$.

Next, we discuss about the smoothness and integrability conditions on the errors. As before, (\ref{hyp:tassio13_eps}) is a regularity condition that enables to use strong approximations for the quantile process. One might be surprised that the moment condition (\ref{hyp:tassio12}) does not show up here, but in fact it is contained in (\ref{uniformmoments}); recall that $\psi_j(\theta_j^*,x)=x$. The integrability condition (\ref{hyp:tassioA3}) is necessary and sufficient for ensuring 
$$
\int_0^1 \frac{B(t)^2}{g^2\{G^{-1}(t)\}}dt<\infty,
$$
from which we see that the limiting random variable in the last claim in Theorem~\ref{prop:tcl} is an almost surely finite random variable. This implies that, as $n \to \infty$,
$$
nW_2^2(G_n,G) \rightsquigarrow  \int_0^1 \frac{B(t)^2}{g^2 \{G^{-1}(t)\}}dt,
$$
with $G_n$ the empirical distribution function on a sample of size $n$ and distribution function $G$. We refer to \cite{delBarrioGineUtzet, SamworthJohnson} for details. Condition (\ref{hyp:tassio12}) is a strong assumption on the tails of $G$ and does not include, e.g.,  {Gaussian} distributions.  {In contrast}, under the less stringent condition
\begin{equation}\label{relaxed}
\int_0^1\int_0^1 \frac{(s\wedge t -st)^2}{g^2\{ G^{-1}(s)\} g^2 \{ G^{-1}(t)\}}dsdt<\infty,
\end{equation}
which is satisfied for  {Gaussian} laws, it can be shown that the limit as $\delta\to 0$
$$
\int_\delta^{1-\delta}\frac{B(t)^2-t(1-t)}{g^2\{G^{-1}(t)\}}dt,
$$
exists in probability and can be expressed as a weighted sum of independent, centered $\chi_1^2$ random variables; see \cite{delBarrioGineUtzet} for details. Then, denoting that kind of limit as
$$
\int_0^{1}\frac{B(t)^2-t(1-t)}{g^2\{G^{-1}(t)\}}dt, 
$$
under some additional tail conditions --- still satisfied by Gaussian distributions; these are conditions (2.10) and (2.22) to (2.24) in the cited reference --- we have that, as $n \to \infty$,
$$
nW_2^2(G_n,G)-c_n \rightsquigarrow  \int_0^1 \frac{B(t)^2-t(1-t)}{g^2 \{ G^{-1}(t) \} }dt,
$$
with 
$$
c_n=\int_{1/n}^{1-1/n} \frac{E B(t)^2}{g^2\{G^{-1}(t)\}}dt. 
$$

A simple look at the proof of Theorem 5.1 shows that under these conditions, instead of (\ref{hyp:tassioA3}), we can conclude that, as $n \to \infty$, 
\begin{equation}\label{centeredCLT}
n A_n(\Theta)-  (J-1) c_n/J^2 \rightsquigarrow  
\frac 1 J\sum_{j=1}^J \int_0^1 \frac{\tilde{B}^2_j(t)- (J-1) t(1-t)/J}{g^2\{ G^{-1}(t)\}  }dt -\frac 1 2 Y^{\top}\Sigma^{-1}Y.
\end{equation}

Our last comment about the assumptions for Theorem~\ref{prop:tcl} concerns the compactness assumption on the parameter space. This may lead in some examples to artificial constraints on the parameter space. However, under some conditions (see, e.g., Corollary~3.2.3 in \cite{vanderVaartWellner})  it is possible to prove that the global minimizer of the empirical criterion lies in a compact neighborhood of the true minimizer. In such cases the conclusion of Theorem~\ref{prop:tcl} would extend  {to} the unconstrained deformation model. 

As a toy example consider the case of deformations by changes in scale, with $J=2$. As above we fix the parameters of, say, the first sample, and consider the family of deformations $\varphi_\sigma(x)=\sigma x$. We assume that the deformation model holds, with the first sample having distribution function $G$ and the second $(\sigma^*)^{-1} G^{-1}$; hence, $\sigma^*$ is the unique minimizer of $U(\sigma)$. We find that 
$$
U_n(\sigma)= \int_0^1 (F_{n,1}^{-1}-\sigma F_{n,2}^{-1})^2/4,
$$
from which we see that almost surely,
$$
\hat{\sigma}_n = \left( \int F_{n,1}^{-1} F_{n,2}^{-1} \right) {\Big /} \left\{ \int  (F_{n,2}^{-1} )^2 \right\} \to \sigma^*
$$ 
and thus the conclusion of Theorem~\ref{prop:tcl} remains valid if we take $\Theta=(0,\infty)$. To avoid further technicalities, we prefer to think of this as a different problem that should be handled in an ad hoc way for each particular example.

Turning back to our goal of assessment of the deformation model (\ref{eq:H0})
based on the observed value of $A_n(\Theta)$,
Theorem~\ref{prop:tcl} gives some insight into the threshold levels for rejection
of the null in the testing problem (\ref{classicalTest}). However, the limiting distribution still 
depends on unknown objects and designing a tractable test requires to estimate the quantiles of this 
distribution. This is the goal of our next result.

We consider bootstrap samples $X^*_{1,j},\ldots,X^*_{m_{n},j}$ 
of  {iid} observations sampled from $\mu^n_j$, write $\mu^*_{m_n,j}$ for the empirical measure on
$X^*_{1,j},\ldots,X^*_{m_n,j}$ and $A_{m_n}^*(\Theta)$ for the minimal alignment cost computed from the bootstrap
samples. We also write $\hat{c}_n(\alpha)$ for the conditional $\alpha$ quantile of $m_n A_{m_n}^*(\Theta)$ given the $X_{i,j}$.

\begin{corollary}\label{cor:bootstrap_th}
Assume that the semiparametric deformation models holds. If $m_n\to \infty$, and $m_n/{n}\to 0$, then under assumptions 
\eqref{hyp:identifTassio}--\eqref{hyp:tassioA3}, we have that
\begin{equation}\label{eq:res_cor_th}
\p \left\{ nA_n( \Theta)> \hat{c}_n(1-\alpha) \right\} \to \alpha. 
\end{equation}
\end{corollary}
Corollary \ref{cor:bootstrap_th} show that the test that rejects $\mathcal{H}_0 : A(\Theta)=0$ (which, as disussed in section 2, is true if and only if the deformation model holds) when $nA_n( \Theta)> \hat{c}_n(1-\alpha)$ is asymptotically of level $\alpha$. It is easy to check that the test is consistent against alternatives that satisfy regularity and integrability assumptions as in Theorem
\ref{prop:tcl}. 

The key to Corollary \ref{cor:bootstrap_th} is that under the assumptions a bootstrap Central Limit Theorem holds for $m_nA_{m_n}^*(\Theta)$. As with Theorem~\ref{prop:tcl}, the integrability conditions on the errors can be relaxed and still have a bootstrap Central Limit Theorem. That would be the case if we replace (\ref{eq:res_cor_th}) by (\ref{relaxed}) and the additional conditions mentioned above under which (\ref{centeredCLT}) holds. Then, the further assumption that the errors have a log-concave distribution and $m_n = O(n^{\rho})$ for some $\rho\in(0,1)$ would be enough to prove a bootstrap Central Limit Theorem, see the comments after the proof of Corollary \ref{cor:bootstrap_th} in the Appendix. In particular, a bootstrap Central Limit Theorem holds for Gaussian tails.
 
\section{Simulations\label{sec:6}}

We present in this section different simulations in order to study  the goodness of fit test we propose in this paper. In this framework, we consider the scale-location family of deformations, i.e., $\theta^*=(\mu^*,\sigma^*)$ and observations such that $X_{i,j}=\mu_j^*+\sigma_j^* \epsilon_{i,j}$ for different distributions of $\epsilon_{i,j}$.

\subsection{Construction of an $\alpha$-level test}

First, we aim at studying the bootstrap procedure which enables to build the test. For this we choose a level $\alpha=0.05$ and aim at estimating the quantile of the asymptotic distribution using a bootstrap method.\vskip .1in
 Let  $B$ be the number of bootstrap samples, we proceed as follows to design a bootstrapped goodness of fit test.
	  		   	\begin{enumerate}
	  		   		\item For all $b \in \{ 1, \ldots,B\}$,
	  		   		\begin{enumerate}
	  		   			\item [1.1.]For $j \in \{ 1, \ldots, J\}$, create a bootstrap sample  $X_{1,j}^{*^{b}}, \ldots, X_{m,j}^{*^{b}}$, with fixed size $m \in \{ 1, \ldots, n\}$ of the first observation sample  $X_{1,j}, \ldots, X_{n,j},$
	  		   			\item [1.2.]Compute  $ (u_{m}^{*b} )^{2}= \displaystyle\inf_{\theta \in \Theta}U_{m}^{*b}(\theta)$. 
	  		   		\end{enumerate}
	  		   		
	  		   	  	\item Sort the values  $ (u_{m}^{*b}   )^{2}$ for $b \in \{ 1,\ldots,B\}$, viz.
	  		   		$ (u_{m}^{*(1)} )^{2} \leqslant \cdots \leqslant  (u_{m}^{*(B)} )^{2}$,
	  		   	    then take $\hat{q}_{m}(1-\alpha)=u_{m}^{*(B(1-\alpha))}$, the quantile of order $1-\alpha$ of the  bootstrap distribution of the statistic $\displaystyle\inf_{\theta \in \Theta}U_{n}(\theta)$. 
	  		   	    \item The test rejects the null hypothesis  if  $nu_{n}^{2} > m [u_{m}^{*\{B(1-\alpha)\}}]^{2}$.
	  		   	\end{enumerate}
	  		   	
Once the test is built, we first ensure that the level of the test has been correctly achieved. For this we repeat the test for large $K$ (here $K=1000$) to estimate the probability of rejection of the test as 
$$
\hat{p}_{n}=\frac{1}{K}\sum_{k=1}^{K} \mathbf{1}{\left[nu_{n,k}^{2} > m\left[u_{m,k}^{*\{B(1-\alpha)\}}\right]^{2}\right]},
$$ 
 {where $\mathbf{1}$ denotes an indicator function.} We present in Table~\ref{tab:frecuenciabajomodelo} these results for different $J$ and several choices for $m=m_n$ depending on the size of the initial sample. 

As expected, the bootstrap method  {makes it possible} to build a test of level $\alpha$ provided the bootstrap sample is large enough. The required size of the sample increases with the number of different distributions $J$ to be tested. 
			\subsection{Power of the test procedure}
			Then we compute the power of previous test for several situations. In particular we must compute the probability of rejection of the null hypothesis under $\mathcal{H}_a$. Hence for several distributions, we test the assumption that the model comes from a warping frame,  when observations from a different distribution called $\gamma$ are observed. The simulations are conducted for the following choices of  {sample size} and for the different distributions:
\begin{itemize}
	  		   		\item [] $J=2$: $\NN (0,1 )$ and $\gamma$;  	  		   		
					\item [] $J=3$ :$\NN (0,1 )$, $\NN (5,2^{2} )$ and $\gamma$; 
	  		   		\item [] $J=5$ : $\NN (0,1 )$, $\NN (5,2^{2} )$, $\NN (3,1 )$, $\NN (1\ldotp5,3^{2} )$ and $\gamma$;
	  		   		\item[] $J=10$: $\NN t(0,1 )$, $\NN (5,2^{2}\ )$, $\NN (3,1\ )$, $\NN (1\ldotp5,3^{2} )$, $\NN (7,4^{2} )$, $\NN (2\ldotp5,0\ldotp5^{2} )$, $\NN (1,1\ldotp5^{2} )$, $\NN (4,3^{2} )$, $\NN (6,5^{2} )$ and $\gamma$;	   	
	  		    \end{itemize}
			    and  also for different choices of  $\gamma$.
			    \begin{itemize}
			    \item [] Exponential distribution with parameter 1, $\mathcal{E}(1)$;
			    \item [] Double exponential (or Laplace) with parameter 1;
			    \item [] Student's $t$ distribution $t_{(3)}$ and $t_{(4)}$ with 3 and 4 degrees of freedom.
			    \end{itemize}
			    
			    All simulations were done for different sample sizes and different bootstrap samples, $n$  and $m_{n}$. The results are presented in Tables~ \ref{tab:potenciaexp}, \ref{tab:potenciadobleexp},  \ref{tab:potenciat3} and \ref{tab:potenciat4}, respectively. 

We observe that the power of the test is very high in most cases. For the Exponential distribution, the power is close to $1$. Indeed  this distribution is very different from the Gaussian distribution since it is  {asymmetric},  {making it easy} to discard the null assumption. The three other distributions do share with the Gaussian the property of symmetry, and yet the power of the test is also close to  {$1$};  {it also increases} with the  {sample size}. Finally, for Student's  {$t$} distribution, the higher the number of degrees of freedom, the more similar it becomes to a Gaussian distribution. This explains why it becomes more difficult for the test to reject the null hypothesis when using a  {Student $t_{(4)}$ rather than a $t_{(3)}$.}

\section*{Appendix} \label{s:append}
\subsection*{A. Proofs of the results in Section~\ref{sec:3}}

\bigskip
\noindent
\textbf{Proof of Theorem~\ref{prop:transportationcost}.}
We set $T_n
={W}_r(\nu_n, \eta)$,  
$T'_n={W}_r(\nu'_n, \eta)$, and let $\Pi_n(\eta)$ be the set of probabilities on $\{1,\ldots,n\}\times \mathbb{R}^d$
with first marginal equal to the discrete uniform distribution on $\{1,\ldots,n\}$ and second marginal equal to
$\eta$. We note that $T_n=\inf_{\pi\in\Pi_n(\eta)} a(\pi)$ if we denote
$$
a(\pi) = \left\{ \int_{\{1,\ldots,n\}\times \mathbb{R}^d} \|Y_i-z\|^rd\pi(i,z)\right\}^{1/r}.
$$
We define similarly $a'(\pi)$ from the $Y'_i$ sample to get $T'_n=\inf_{\pi\in\Pi_n(\eta)} a'(\pi)$. But
then, using the inequality $ \left\vert\,  \left\Vert a \right\Vert -\left\Vert b \right\Vert\, \right\vert \leq \left\Vert a-b \right\Vert $,  {we get}
\begin{equation*}
|a(\pi)-a'(\pi)|\leq  \left\{ \int_{\{1,\ldots,n\}\times \mathbb{R}^d} \|Y_i-Y_i'\|^rd\pi(i,z)\right\}^{1/r}=
\left( \frac 1 n  \sum_{i=1}^n\|Y_i-Y_i'\|^r\right)^{1/r}. \label{eq:boots1}
\end{equation*}
This implies that
\begin{equation*}|T_n-T_n'|^r\leq  \frac 1 n  \sum_{i=1}^n\|Y_i-Y_i'\|^r.\label{eq:boots2}
\end{equation*}

If we take now $(Y,Y')$ to be an optimal coupling of $\nu$ and $\nu'$, so that ${\rm E} ( \|Y-Y'\|^r ) = {W}_r^r(\nu,\nu')$
and $(Y_1,Y'_1),\ldots,(Y_n,Y'_n)$ to be  {iid} copies of $(Y,Y')$, we see that
for the corresponding realizations of $T_n$ and $T_n'$, we have
$$
\E ( |T_n-T_n'|^r) \leq  \frac 1 n  \sum_{i=1}^n \E\left(\|Y_i-Y_i'\|^r\right)={W}_r(\nu,\nu')^r.
$$
But this shows that ${W}_r\{\mathcal{L}(T_n),\mathcal{L}(T_n')\}\leq {W}_r(\nu,\nu')$, as claimed.
\hfill $\Box$

\bigskip
\noindent
\textbf{Proof of Theorem~\ref{prop:rvariation}.}
We write $V_{r,n}=V_r(\nu_{n_1,1},\ldots,\nu_{n_J,J})$ and $V'_{r,n}=V_r(\nu'_{n_1,1},\ldots,\nu'_{n_J,J})$.
We note that
$$V_{r,n}^r=\inf_{\pi\in \Pi(U_1,\ldots,U_J)} \int T(i_1,\ldots,i_J) d\pi(i_1,\ldots,i_J),$$
where $U_j$ is the discrete uniform distribution on $\{1,\ldots,n_j\}$ and 
$$
T(i_1,\ldots,i_J)=\min_{z\in \mathbb{R}^d} 
\frac 1 J \sum_{j=1}^J\|Y_{i_j,j}-z\|^r. 
$$
We write $T'(i_1,\ldots,i_J)$ for the equivalent function computed from the
$Y'_{i,j}$s. Hence we have
$$|T(i_1,\ldots,i_J)^{1/r}-T'(i_1,\ldots,i_J)^{1/r}|^r\leq\frac 1 J \sum_{j=1}^J \|Y_{i_j,j}-Y'_{i_j,j}\|^r,$$
which implies
\begin{eqnarray*}
\lefteqn{\left|\left\{\int T(i_1,\ldots,i_J) d\pi(i_1,\ldots,i_J)  \right\}^{1/r}- \left\{ \int T(i_1,\ldots,i_J) 
d\pi(i_1,\ldots,i_J)  \right\}^{1/r}\right|^r}\\
\quad &\leq&  \int\frac 1 J \sum_{j=1}^J \|Y_{i_j,j}-Y'_{i_j,j}\|^r d\pi(i_1,\ldots,i_J)\\
\quad &=&\frac 1 J \sum_{j=1}^J\int \| Y_{i_j,j}-Y'_{i_j,j}\|^rd\pi(i_1,\ldots,i_J)
=\frac 1 J \sum_{j=1}^J\left( \frac 1 {n_j} \sum_{i=1}^{n_j} \| Y_{i,j}-Y'_{i,j}\|^r\right).
\end{eqnarray*}
 {Therefore,}
$$
|V_{r,n}-V'_{r,n}|^r\leq 
\frac 1 J \sum_{j=1}^J\left( \frac 1 {n_j} \sum_{i=1}^{n_j} \| Y_{i,j}-Y'_{i,j}\|^r\right).
$$
If we take $(Y_j,Y'_j)$ to be an optimal coupling of $\nu_j$ and $\nu_j'$ and $(Y_{1,j},Y'_{1,j}),\ldots,$
$(Y_{n_j,j},Y'_{n_j,j})$ to be iid copies of $(Y_j,Y'_j)$ for all $j \in \{ 1,\ldots,J\}$, then
we obtain 
$$
\E\left(|V_{r,n}-V'_{r,n}|^r \right) \leq 
\frac 1 J \sum_{j=1}^J\left\{ \frac 1 {n_j} \sum_{i=1}^{n_j} \E\left(\| Y_{i,j}-Y'_{i,j}\|^r\right)\right\} = \frac 1 J \sum_{j=1}^J {W}_{r}^r(\nu_j,\nu'_j).
$$
The conclusion follows.
\hfill $\Box$

\bigskip
\noindent
\textbf{Proof of Theorem~\ref{th:bootstrap}.}
We argue as in the proof of Theorem~\ref{prop:rvariation} and write
$$
A_{n,r}(\mathcal{G})=\inf_{\varphi \in\mathcal{G}} \left\{ \inf_{ \pi\in \Pi(U_1,\ldots,U_J)}  
\int T(\varphi ;i_1,\ldots,i_J) d\pi(i_1,\ldots,i_J)\right\},
$$
where 
$$
T(\varphi  ;i_1,\ldots,i_J)=\min_{y\in\mathbb{R}} \frac 1 J \sum_{j=1}^J \|Z_{i_j,j}(\varphi_j)-y\|^r.
$$
We write $T'(\varphi  ;i_1,$ $\ldots,i_J)$ for the same function computed on the $Z'_{i,j}(\varphi_j)$'s.
Now, from the  fact $\|Z_{i,j}(\varphi_j)-Z'_{i,j}(\varphi_j)\|^r\leq L^r \|X_{i,j}-X'_{i,j}\|^r$ 
we see that
$$
|T(\varphi ;i_1,\ldots,i_J)^{1/r}-T'(\varphi ;i_1,\ldots,i_J)^{1/r}|^r\leq 
L^r\frac 1 J \sum_{j=1}^J \|X_{i_j,j}-X'_{i_j,j}\|^r
$$
and, as a consequence,
$$
|V_r\{\mu_{n,1}(\varphi_1),\ldots,\mu_{n,J}(\varphi_J)\}-V_r \{\mu'_{n,1}(\varphi_1),\ldots,\mu'_{n,J}(\varphi_J)\}|^r
\leq \frac {L^r} J \sum_{j =1}^J \sum_ {i_j=1}^{n_j} \frac{1}{n_j} \|X_{i_j,j}-X'_{i_j,j}\|^r
$$
which implies 
$$
|\{A_{n,r}(\mathcal{G})\}^{1/r} - \{A'_{n,r}(\mathcal{G})\}^{1/r}|^r\leq  \frac {L^r} J \sum_{j=1}^J 
\left({ \frac 1 {n_j}\sum_{i=1}^{n_j}\|X_{i,j}-X'_{i,j}\|^r}
\right).$$
If, as in the proof of Theorem~\ref{prop:rvariation}, we assume that $(X_{i,j},X'_{i,j})$ with $i \in \{1,\ldots,n_j\}$ are iid copies of an optimal coupling for $\mu_j$ and $\mu_j'$, with different samples independent from each other we obtain that
$$
\E\left[ |\{A_{n,r}(\mathcal{G})\}^{1/r} - \{A'_{n,r}(\mathcal{G})\}^{1/r}|^r\right]\leq \frac {L^r} J \sum_{j=1}^J {W}_r^r(\mu_j,\mu_j').
$$
 {This concludes the argument.}
\hfill $\Box$

\subsection{Proofs for results from Sections~\ref{sec:4} and \ref{sec:5}}

We provide here proofs of the main results in Sections~\ref{sec:4} and \ref{sec:5}. Our approach
relies on the consideration the processes defined, for all $\varphi\in\mathcal{G}$, by
$$
C_n(\varphi)=\sqrt{n} \, \{U_n(\varphi)-U(\varphi)\} \quad \mbox{and} \quad C(\varphi)=\frac 1 J \sum_{j=1}^J c_{j}(\varphi),
$$
where $U_n(\varphi)=V_2^2\{\mu_{n,1}(\varphi_1), \dots, \mu_{n,J}(\varphi_J)\}$, 
$U(\varphi)=V_2^2\{ \mu_1(\varphi_1), \dots, \mu_J(\varphi_J)\}$,
$$c_{j}(\varphi) =2 \int_0^1  \varphi_j^\prime \circ F_j^{-1} \{\varphi_j\circ F_j^{-1}-F_B^{-1}(\varphi)\} \frac{B_{j}}{f_j\circ F_j^{-1}}$$
and $B_1, \ldots, B_J$ are independent standard Brownian bridges on $(0,1)$. We prove below that the empirical
deformation cost process $C_n$ converges weakly to $C$ as random elements in $L^\infty(\mathcal{G})$, 
the space of bounded,  {real-valued} functions on $\mathcal{G}$. 
Theorem~\ref{th:testgen} will follow as a corollary of this result.

We will make frequent use in this section of the following technical  {lemma}, 
which follows easily from the triangle and  {H\"older's} inequalities. We omit the proof.
\begin{lem}
\label{lem:ln}
Under Assumption \eqref{hyp:tassio12theta},
\begin{enumerate}
\item[(i)] 
$$
\sup_ {\varphi_j \in \mathcal{G}_j} \sqrt{n} \int_0^{1/n} \left(\varphi_j\circ F_j^{-1} \right)^2 \rightarrow 0,\quad  
\sup_ {\varphi_j \in \mathcal{G}_j} \sqrt{n} \int_{1-1/n}^1 \left(\varphi_j\circ F_j^{-1} \right)^2  \rightarrow 0.
$$
\item[(ii)] 
$$
\sup_ {\varphi_j \in \mathcal{G}_j} \sqrt{n} \int_0^{1/n} \left(\varphi_j\circ F_{n,j}^{-1} \right)^2 \rightarrow 0, \quad  
\sup_ {\varphi_j \in \mathcal{G}_j} \sqrt{n} \int_{1-1/n}^1 \left(\varphi_j\circ F_{n,j}^{-1}\right)^2  \rightarrow 0
$$ 
in probability.
\item[(iii)] If moreover \eqref{hyp:Rajput} holds, then for all $j, k \in \{ 1, \ldots, J\}$,
$$
 \int_0^1 \frac{\sqrt{t(1-t)}}{f_k\{ F_k^{-1} (t) \}} \sup_{\varphi_j \in \mathcal{G}_j}\Big\vert  \varphi_j \{ F_j^{-1}( t ) \}\Big\vert dt < \infty.
$$
\end{enumerate}
\end{lem}

\begin{theo}
Under assumptions \eqref{hyp:tassio13}--\eqref{hyp:tassio12theta}, $C_n$ and $C$ have almost surely trajectories in 
$L^\infty(\mathcal{G})$. Furthermore, $C$ is a tight Gaussian random elemnt and $C_n$ converges weakly to $C$ in $L^\infty(\mathcal{G})$.
\end{theo}
\noindent 
\textit{Proof.} We start  {by} noting that 
$$
U_n(\varphi)=\frac 1 J \sum_{j=1}^J \int_0^1 \left\{\varphi_j \circ F_{n,j}^{-1}-F_{n,B}^{-1}(\varphi) \right\}^{2},\quad
U(\varphi)=\frac 1 J \sum_{j=1}^J \int_0^1 \left\{\varphi_j \circ F_{j}^{-1}-F_B^{-1}(\varphi) \right\}^{2}
$$ with 
$$
F_{n,B}^{-1}(\varphi)=\frac 1 J\sum_{j=1}^J \varphi_j\circ F_{n,j}^{-1},\quad 
F_B^{-1}(\varphi)=\frac 1 J\sum_{j=1}^J \varphi_j\circ F_j^{-1}.
$$
Now,~(\ref{hyp:tassio12theta})
implies that 
$$
\sup_{\varphi_j\in\mathcal{G}_j} \int_0^1 \left(\varphi_j\circ F_j^{-1}\right)^2 < \infty.
$$
Similarly, assumption (\ref{hyp:regul}) implies $K_j = \sup_{\varphi_j\in\mathcal{G}_j, x\in (c_j,d_j)}|\varphi_j'(x)|<\infty$. 
Noting that 
$$
\int_0^1 \left(\varphi_j\circ F_{n,j}^{-1} \right)^2\leq 2 \int_0^1 \left(\varphi_j\circ F_{j}^{-1}\right)^2 + 2 K_j^2 \int_0^1 \left(F_{n,j}^{-1} - F_{j}^{-1} \right)^2,
$$  
we see that 
$$
\sup_{\varphi_j\in\mathcal{G}_j} \int_0^1 \left(\varphi_j\circ F_{n,j}^{-1}\right)^2<\infty,\; \mbox{a.s.}$$ and, with little additional effort,
conclude that $C_n$ has almost surely bounded trajectories. Furthermore, writing 
$$
d_{j,k}(\varphi) = \int_0^1  \varphi_j^\prime \circ F_j^{-1} \frac{B_{j}}{f_j\circ F_j^{-1}} \varphi_k\circ F_k^{-1},
$$
we see that 
for  $\varphi, \rho \in \mathcal{G}$
\begin{align*}
|d_{j,k}(\varphi) &- d_{j,k}(\rho)| \\
&\leq  \|\varphi_j'-\rho_j'\|_\infty 
\left|\int_0^1 \frac{B_{k}}{f_k\circ F_k^{-1}} \varphi_k\circ F_k^{-1} \right|  
+ \left|\int_0^1 \rho_j'\circ F_j^{-1} \frac{B_{k}}{f_k\circ F_k^{-1}} \left(\varphi_k\circ F_k^{-1}-\rho_k\circ F_k^{-1} \right) \right|\\
&\leq  \|\varphi_j'-\rho_j'\|_\infty  \sup_{\varphi_k \in \mathcal{G}_k}\left|\int_0^1 \frac{B_{k}}{f_k\circ F_k^{-1}} \varphi_k\circ F_k^{-1} \right| +\sup_{(c_j,d_j)} |\rho_j^\prime |  \left({\int_0^1 \Big| \frac{B_{k}}{f_k\circ F_k^{-1}}}\Big|^{q} \right)^{1/q}
\left({\int_0^1 \left|\varphi_k\circ F_k^{-1}-\rho_k\circ F_k^{-1}\right|^{p_0}}\right)^{1/p_0}
\end{align*}
But using (iii) of Lemma~\ref{lem:ln}, we deduce that
$$
\E  \left( \sup_{\varphi_k \in \mathcal{G}_k}\bigg|\int_0^1 \frac{B_{k}}{f_k\circ F_k^{-1}} \varphi_k\circ F_k^{-1}\bigg| \right)
\leq \int_0^1 \frac{\sqrt{t(1-t)}} {f_k\{F_k^{-1}(t)\}} \sup_{\varphi_j  \in \mathcal{G}_j}  \vert \varphi_j \{F_j^{-1}( t)\}\vert dt < \infty.
$$
Hence, almost surely, 
$$
\sup_{\varphi\in\mathcal{G}}\bigg|\int_0^1 \frac{B_{j}}{f_j\circ F_j^{-1}} \varphi_j\circ F_j^{-1} \bigg|<\infty.
$$
Furthermore, from assumption~(\ref{hyp:Rajput}), we get that, almost surely,
$$
\int_0^1 \left({\frac{B_{j}}{f_j\circ F_j^{-1}}}\right)^{q}< \infty, 
$$
and thus, for some almost surely finite random variable $T$, $\vert d_{j,k} ( \varphi) - d_{j,k} (\rho) \vert
\leq T \Vert \varphi - \rho \Vert_{\mathcal{G}}$ for  $\varphi$, $\rho \in \mathcal{G}$.
From this conclude that the trajectories of $C$ are a.s. bounded, uniformly continuous 
functions on  $\mathcal{G}$, endowed with the norm  $\left\Vert \cdot \right\Vert_{\mathcal{G}}$
introduced in (\ref{hyp:regul}). In particular, $C$ is a tight random element in ${L}^\infty(\mathcal{G})$; see, e.g.,  {pp.} 39--41 in~\cite{vanderVaartWellner}.

From this point  {on,} we pay attention to the quantile processes defined, for all $j \in \{ 1, \ldots, J\}$ and $t \in (0,1)$, by
$$
\rho_{n,j}(t)=\sqrt{n}\,f_j \left\{F_j^{-1}(t)\right\} \left\{F_{n,j}^{-1}(t)-F_{j}^{-1}(t)\right\}.
$$
A trivial adaptation of Theorem~2.1 on p.~381 of~\cite{csorgohorvath93} shows that,
under (\ref{hyp:tassio13}), there exist, on a rich enough probability space, independent versions of
$\rho_{n,j}$ and independent families of Brownian bridges $\{B_{n,j}\}_{n=1}^\infty$ for all $j \in \{1,\ldots,J\}$,
satisfying 
\begin{equation}\label{csorgohorvath}\tag{A.1}
n^{1/2 - \nu} \sup_{1/n \leq t \leq 1- 1/n} \frac{| \rho_{n,j}(t) - B_{n,j}(t) 
|}{\left\{ t (1-t) \right\} ^\nu}
= \begin{cases} 
O_p(\ln n)& \mbox{if } \nu=0,\\
O_p(1)& \mbox{if } 0< \nu \leq 1/2.
\end{cases}
\end{equation}
We work, without loss of generality, with these versions of $\rho_{n,j}$ and $B_{n,j}$. We show now that
\begin{equation}\label{eq:step2fin}\tag{A.2}
\sup_{\varphi  \in \mathcal{G}} | C_n(\varphi)- \hat{C}_n(\varphi) |  \rightarrow 0 \text{ in probability}
\end{equation}
with 
$$
\hat{C}_n(\varphi) = \frac{1}{J} \sum_{j=1}^J c_{n,j}\left(\varphi \right) \quad \mbox{and} \quad
c_{n,j}(\varphi) = 2 \int_0^1  \varphi_j^\prime \circ F_j^{-1} \left\{ \varphi_j\circ F_j^{-1}-F_B^{-1}(\varphi)\right\}  \frac{B_{n,j}}{f_j\circ F_j^{-1}}.
$$
To check this, we note that some simple algebra yields
$$
C_n(\varphi)=\frac{2}{J}\sum_{j=1}^J \tilde{c}_{n,j} + \frac{1}{J} \sum_{j=1}^J\tilde{r}_{n,j} 
$$ 
with 
$$
\tilde{c}_{n,j} = \sqrt{n}\int_0^1 \left(\varphi_j\circ F_{n,j}^{-1}-\varphi_j\circ F_j^{-1} \right)
\left( \varphi_j\circ F_j^{-1}-F_B^{-1}(\varphi)\right),
$$
$$
\tilde{r}_{n,j} = \sqrt{n}\int_0^1 \left[ \left(\varphi_j\circ F_{n,j}^{-1}-\varphi_j\circ F_j^{-1}\right) - \left\{F_{n,B}^{-1}(\varphi)-F_B^{-1}(\varphi) \right\} \right]^2.
$$
From the elementary inequality $(a_1+\cdots+a_J)^2\leq Ja_1^2+\cdots+Ja_J^2$, we get that
$$
\frac{1}{J} \sum_{j=1}^J\tilde{r}_{n,j}\leq \frac{4\sqrt{n}}{J} \sum_{j=1}^J \int_0^1 \left(\varphi_j\circ F_{n,j}^{-1}-\varphi_j\circ F_{j}^{-1} \right)^2\leq 
\frac{4\sqrt{n}}{J} \sum_{j=1}^J K_j \int_0^1 \left(F_{n,j}^{-1}-F_{j}^{-1} \right)^2,
$$
with $K_j = \sup_{\varphi_j\in\mathcal{G}_j, x\in (c_j,d_j)}|\varphi_j'(x)|<\infty$, as above.
Now we can use (\ref{hyp:tassio12}) and argue as in the proof of Theorem 2 in \cite{MR2435470} to conclude
that $\sqrt{n}\int_0^1 (F_{n,j}^{-1}-F_{j}^{-1})^2\to 0$ in probability and, as a consequence, that
\begin{equation}
 \label{eq:intermediate}\tag{A.3}
 \sup_{\varphi  \in \mathcal{G}} \bigg| C_n(\varphi)- 
 \frac{1}{J} \sum_{j=1}^J \tilde{c}_{n,j}\left(\varphi \right) \bigg|  \rightarrow 0 \mbox{ in probability}.
\end{equation}
Furthermore, the Cauchy--Schwarz inequality  {implies} that
\begin{multline*}
n\left[\int_0^{ {1/n} } \left(\varphi_j\circ F_{n,j}^{-1}-\varphi_j\circ F_j^{-1} \right)\left\{\varphi_j\circ F_j^{-1}-F_B^{-1}(\varphi) \right\} \right]^2\\
\leq \sqrt{n} \int_0^{ {1/n} } \left(\varphi_j\circ F_{n,j}^{-1}-\varphi_j\circ F_j^{-1} \right)^2
 \sqrt{n} \int_0^{ {1/n} } \left\{\varphi_j\circ F_j^{-1}-F_B^{-1}(\varphi)\right\}^2
\end{multline*}
and using (i) and (ii) of Lemma~\ref{lem:ln}, the two factors converge to zero uniformly in $\varphi$.
A similar argument works for the upper tail and allows to conclude that
we can replace $\tilde{c}_{n,j}(\varphi)$ in~(\ref{eq:intermediate}) with 
$$
\tilde{\tilde{c}}_{n,j}(\varphi)=
2 \sqrt{n}\int_{ {1/n} }^{1- {1/n} } \left(\varphi_j\circ F_{n,j}^{-1}-\varphi_j\circ F_j^{-1}\right) \left\{ \varphi_j\circ F_j^{-1}-F_B^{-1}(\varphi)\right\}.
$$
Moreover,
$$
\sup_{\varphi  \in \mathcal{G}}  \bigg| \int_0^{1/n}\varphi_j^\prime\circ F_j^{-1}\frac{B_{n,j} }{ {f_j\circ F_j^{-1}}} 
\left\{\varphi_j \circ  F_j^{-1}-F^{-1}_B(\varphi)\right\} \bigg|
\leq 
K_j \int_0^{1/n} \Big| \frac{B_{n,j}}
{f_j\circ F_j^{-1}}\Big| \sup_{\varphi  \in \mathcal{G}}  \left| \left\{\varphi_j \circ  F_j^{-1}-F^{-1}_B(\varphi)\right\} \right| 
$$
and by (iii) of Lemma~\ref{lem:ln} and Cauchy--Schwarz's inequality,
$$
\E  \left[ \int_0^\frac{1}{n} \Big| \frac{B_{n,j}}
{f_j\circ F_j^{-1}}\Big| \sup_{\varphi  \in \mathcal{G}}  \left| \left\{\varphi_j \circ  F_j^{-1}-F^{-1}_B(\varphi)\right\} \right|  \right]
 \leq  
\int_0^\frac{1}{n} \frac{ \sqrt{t(1-t)}  }{ {f_j\{F_j^{-1}(t)\}}}\sup_{\varphi  \in \mathcal{G}}  \left| \varphi_j 
\{ F_j^{-1}( t )\} - F^{-1}_B(\varphi)
( t) \right| dt   \rightarrow 0.
$$
Hence,  
$$
\sup_{\varphi  \in \mathcal{G}}  \left| \int_0^\frac{1}{n}\varphi_j^\prime\circ F_j^{-1}\frac{B_{n,j} }{ {f_j\circ F_j^{-1}}} \left\{\varphi_j \circ  F_j^{-1}-F^{-1}_B(\varphi)\right\} \right|\to 0
$$
in probability and similarly for the right tail. 
Now,  for every $t\in (0,1)$ we have 
\begin{equation}
\label{taylor}\tag{A.4}
\varphi_j\circ F_{n,j}^{-1}(t)- \varphi_j\circ F_{j}^{-1}(t)= 
\varphi_j^\prime \{ K_{n,\varphi_j}(t)\} \{F_{n,j}^{-1}(t)- F_{j}^{-1}(t)\}
\end{equation}
for some $K_{n,\varphi_j}(t)$ between $F_{n,j}^{-1}(t)$ and $F^{-1}(t)$. Therefore, recall~\eqref{taylor}, to prove~(\ref{eq:step2fin}) it suffices to show that
\begin{multline}\label{finalapp}\tag{A.5}
\sup_{\varphi\in\mathcal{G}}\bigg|  
\int_{1/n}^{1- 1/n} \varphi_j^\prime \{F_j^{-1}(t)\}\frac{B_{n,j}(t) }{ {f_j\{F_j^{-1}(t)\}}} 
\left[\varphi_j \{F_j^{-1}(t)\}-F^{-1}_B(\varphi)(t)\right] dt \\
- \int_{1/n}^{1- {1/2} }\varphi_j^\prime \{K_{n,\varphi_j}(t)\} \frac{\rho_{n,j}(t) }{ {f_j \{F_j^{-1}(t)\}}} \left[\varphi_j \{F_j^{-1}(t)\}-F^{-1}_B(\varphi)(t)\right] dt  \bigg|\to 0
\end{multline}
in probability. To check it, we take $\nu\in(0,1/2)$ 
in (\ref{csorgohorvath}) to get
\begin{multline}\label{finalapp1}\tag{A.6}
\int_{1/n}^{1-{1/n} } \frac{|\rho_{n,j}(t)- B_{n,j}(t)|}{f_j\{F_j^{-1}(t)\}}
\sup_{\varphi  \in\mathcal{G}} \big| \varphi_j \{F_j^{-1}( t )\} -F^{-1}_B
( \varphi ) (t ) \big| dt \\  
\leq n^{\nu -{1}/{2}  }O_P(1)
\int_{1/n}^{1-{1/n} } \frac{ \{t( 1-t)\} ^\nu}{f_k \{F_k^{-1}(t)\}}\sup_{\varphi  \in\mathcal{G}} \big| 
\varphi_j \{F_j^{-1}( t )\} -F^{-1}_B( \varphi )( t) \big|  dt\to 0
\end{multline}
in probability using dominated convergence and (iii) of Lemma~\ref{lem:ln}.
We observe next that, for each $t\in (0,1)$, 
$$
\sup_{ \varphi_j \in  \mathcal{G}_j}| K_{n,\varphi_j}(t) - F_j^{-1}(t)| \rightarrow 0
$$
a.s., since
$ K_{n,\varphi_j}(t)$ lies between $F_{n,j}^{-1}(t)$ and $F_{j}^{-1}(t)$. Therefore, using (\ref{hyp:regul}) 
we see that 
$$
\sup_{ \varphi_j\in  \mathcal{G}_j }| \varphi_j^\prime \{K_{n,\varphi_j}(t)\} - \varphi_j^\prime \{F_j^{-1}(t)\}| \rightarrow 0
$$ 
a.s. while, on the other hand, $\sup_{ \varphi_j\in  \mathcal{G}_j }| \varphi_j^\prime \{K_{n,\varphi_j}(t)\} - \varphi_j^\prime \{F_j^{-1}(t)\}|\leq 2K_j$.
But then, by dominated convergence we get that
$$ \E \Big[\sup_{ \varphi_j\in  \mathcal{G}_j } \left| \varphi_j^\prime \{K_{n,\varphi_j}(t)\} - \varphi_j^\prime \{F_j^{-1}(t)\} \right|^2 \Big]  \rightarrow 0.
$$
Since by (iii) of Lemma~\ref{lem:ln} we have that $t\mapsto {\sqrt{t(1-t)}}/[{f_j\{F_j^{-1}(t)\}}] \sup_{\varphi\in\mathcal{G}}|\varphi_j\{F_j^{-1}(t)\}-
F_B^{-1}(\varphi)(t) |$ is integrable we conclude that 
$$
\E \sup_{\varphi\in \mathcal{G}} 
\int_{1 /n}^{1- {1/n} } \big| \varphi_j^\prime \{K_{n,\varphi_j}(t)\} - \varphi_j^\prime \{F_j^{-1}(t)\} \big|
\frac {|B_{n,j}(t)|}{f_j\{F_j^{-1}(t)\}} \big|\varphi_j\{F_j^{-1}(t)\}-
F_B^{-1}(\varphi)(t) \big|dt
$$
tends to 0 as $n\to\infty$ and, consequently,
$$
\sup_{\varphi\in \mathcal{G}} \int_{1 /n}^{1- {1/n} } \big| \varphi_j^\prime \{K_{n,\varphi_j}(t)\} - \varphi_j^\prime \{F_j^{-1}(t)\} \big|
\frac {|B_{n,j}(t)|}{f_j\{F_j^{-1}(t)\}} \big|\varphi_j\{F_j^{-1}(t)\}-
F_B^{-1}(\varphi)(t) \big|dt$$
vanishes in probability. Combining this fact with (\ref{finalapp1}) we prove (\ref{finalapp}) and, as a consequence, (\ref{eq:step2fin}). Finally, observe that for every integer $n\geq 1$, $C$ has the same law as $\hat{C}_n$. This completes the proof. 
\hfill $\Box$

\bigskip
\noindent
{\textbf{Proof of Theorem~\ref{th:testgen}.}
From Skohorod Theorem; see, e.g., Theorem 1.10.4 in~\cite{vanderVaartWellner}, we know that there exists  on 
some probability space versions of $C_n$ and $C$ for which convergence of $C_n$ to $C$ holds almost surely. 
From now on, we place us on this space and observe that
\begin{equation}
\label{ineq:bornesup}\tag{A.7}
\sqrt{n} \, \{A_n(\mathcal{G})-A(\mathcal{G})\}\leq \sqrt{n} 
\inf_\Gamma U_n - \sqrt{n} \inf_\Gamma  U = \inf_{\varphi\in\Gamma} C_n(\varphi).
\end{equation}
Furthermore, if we consider the (a.s.) compact set $\Gamma_n = \lbrace  \varphi \in {\mathcal{G}}:\, U\left( \varphi  \right) \leq \inf_{{\mathcal{G}}} U + {2 \left\Vert C_n \right\Vert_\infty}/{\sqrt{n}} \rbrace$, then, if $\varphi \notin \Gamma_n$, $U_n\left( \varphi \right)\geq \inf_{\mathcal{G}} U + \left\Vert  C_n \right\Vert_\infty/{\sqrt{n}}$, while if $\varphi \in \Gamma$, then, $U_n\left( \varphi \right) \leq \inf_{\mathcal{G}} U + {\left\Vert C_n \right\Vert_\infty}/{\sqrt{n}}$. Thus, necessarily, $ \inf_ { \mathcal{G}} U_n=\inf_ { \Gamma_n} U_n=\inf_ { \Gamma_n} (U_n-U+U)\geq \inf_ { \Gamma_n} (U_n-U)+\inf_{\Gamma_n} U = \inf_ { \Gamma_n} (U_n-U)+\inf_{\Gamma} U$. Together with \eqref{ineq:bornesup}, this entails
\begin{equation*}
\label{eq:encadrement}
\inf_{\varphi\in\Gamma_n} C_n(\varphi)\leq\sqrt{n}\, \{A_n(\mathcal{G})-A(\mathcal{G})\}\leq  \inf_{\varphi\in\Gamma} C_n(\varphi)
\end{equation*}
Note that for the versions that we are considering $\|C_n-C\|_\infty\to 0$~almost surely. In particular, this implies
that $\inf_{\Gamma} C_n\to \inf_{\Gamma} C$~almost surely. Hence, the proof will be complete if we 
show that almost surely,
\begin{equation}\label{finalpart}\tag{A.8}
\inf_{\Gamma_n} C_n\to \inf_{\Gamma} C.
\end{equation}
To check this last point, consider a sequence $\varphi_n\in\Gamma_n$ such that 
$C_n(\varphi_n)\leq \inf_{\Gamma_n}C_n+ {1/n} $.
By compactness of $\mathcal{G}$, taking subsequences if necessary, $\varphi_n\to \varphi_0$ for some $\varphi_0\in\mathcal{G}$. Continuity of $U$ yields $U(\varphi_n)\to U(\varphi_0)$ and as a consequence, that $U(\varphi_0)\leq \inf_{\mathcal{G}} U$, i.e., $\varphi_0\in\Gamma$ almost surely, Furthermore,   
$$
\big| C_n(\varphi_n) - C(\varphi_0) \big|
\leq \left\Vert C_n-C  \right\Vert_\infty + 
\left\vert C\left(\varphi_n\right)- C \left( \varphi_0 \right)\right\vert \to 0.
$$
This shows that
\begin{equation*} \label{eq:b_sup}
\liminf \inf_{\Gamma_n} C_n \geq C\left( \varphi_0 \right)  \geqslant  \inf_\Gamma C
\end{equation*}
and yields (\ref{finalpart}). This completes the proof.
\hfill $\Box$

\bigskip
\medskip
\noindent
\textbf{Proof of Corollary~\ref{cor:bootstrap}.}
In Theorem~\ref{th:bootstrap}, take $\mu'_j = \mu_{n,j} $. 
Then, writing $\mathcal{L}^*$ for the conditional law given the $X_{i,j}$, the result of Theorem~\ref{th:bootstrap} reads 
\begin{equation*} 
{W}_2^2\left[\mathcal{L}\big[\{A_{m_n}(\mathcal{G})\}^{1/2}\big],\mathcal{L}^*\big[\{A^*_{m_n}(\mathcal{G})\}^{1/2}\big]\right]
\leq L^2 \frac 1 J \sum_{j=1}^J {W}_2^2(\mu_j,\mu_{n,j}),
\end{equation*}
with $L=\sup_{\varphi \in \mathcal{G} }\Vert \varphi^\prime_j\Vert_{\infty}<\infty$. Since ${W}_r\{\mathcal{L}(aX+b),\mathcal{L}(aY+b)\}=a{W}_r\{\mathcal{L}(X),\mathcal{L}(Y)\}$ for $a>0,$ $b\in\mathbb{R}$, the  {latter} bound gives
$$ 
{W}_2^2\left[\mathcal{L}\big[\sqrt{m_n}\big[ \{A_{m_n}(\mathcal{G})\}^{1/2}- \{A(\mathcal{G})\}^{1/2}\big]\big],
\mathcal{L}^*\big[\sqrt{m_n}\big[ \{A^*_{m_n}(\mathcal{G})\}^{1/2}-\{A(\mathcal{G})\}^{1/2}\big]\big]\right]
\leq L^2 \frac{{m_n}}{\sqrt{n}} \frac 1 J \sum_{j=1}^J  \sqrt{n}\,{W}_2^2(\mu_j,\mu_{n,j}).
$$
As noted in the proof of Theorem~\ref{th:testgen}, the assumptions imply that $\sqrt{n}\,{W}_2^2(\mu_j,\mu_{n,j})$ vanishes in probability. Also, Theorem~\ref{th:testgen} and the delta method yield that 
$$
\sqrt{m_n}\Big[ \{A_{m_n}(\mathcal{G})\}^{1/2}-\{A(\mathcal{G})\}^{1/2}\Big] \rightsquigarrow \frac{1}{2\{A(\mathcal{G})\}^{1/2}}\gamma,$$
with $\gamma$ the limiting law there, which, combined to the above bound, shows that 
$$
\sqrt{m_n} \, \Big[ \{A^*_{m_n}(\mathcal{G})\}^{1/2}-\{A(\mathcal{G})\}^{1/2}\Big] \rightsquigarrow \frac{1}{2\{A(\mathcal{G})\}^{1/2}}\gamma
$$
in probability. A further use of the delta method yields
$$
\sqrt{m_n} \, \Big\{ A^*_{m_n}(\mathcal{G})-A(\mathcal{G})\Big\} \rightsquigarrow \gamma
$$
in probability. The result follows now from Lemma 1 in \cite{JanssenPauls}.
\hfill $\Box$

\bigskip
\noindent
\textbf{Proof of Theorem~\ref{prop:tcl}.} 
We assume for simplicity that $p=1$. The general case follows with straightforward changes. 
Let us observe that
$$U_n(\theta)=\frac 1 J\sum_{j=1}\int_0^1 \bigg\{\psi_j(\theta_j,G_{n,j}^{-1})-{ \frac 1 J  \sum_{k=1}^J} \psi_k(\theta_k,G_{n,k}^{-1}) \bigg\}^2, $$
with $G_{n,j}$ the empirical distribution function on the $\varepsilon_{i,j}$s, which are  {iid} $G$. A similar expression, replacing
$G_{n,j}$ with $G$ is valid for~$U(\theta)$. Then (\ref{hyp:regul_th}) implies that $\sup_{\theta} |U_n(\theta)-U(\theta)|\to 0$, from which, recall~\eqref{hyp:identifTassio}, it follows that $\hat{\theta}_n\to \theta^*$ in probability. Note that the second part in  Assumption~(\ref{hyp:regul_th}) is a technical  {condition which} ensures that, when considering a Taylor expansion in the integral of $U_n(\theta)$, the remainder term in $\psi^{'}_j(\lambda,H_{n,j}^{-1})-\psi^{'}_j(\lambda,G_{j}^{-1})$ for any $H_{n,j}^{-1}$ lying between $G_{n,j}^{-1}$ and $G_j^{-1}$ (obtained through a Taylor expansion) goes uniformly to zero. 

From (\ref{jointsmoothness}) we have that $U_n$ is a $C^2$ function
whose derivatives can be computed by differentiation under the integral sign. This implies that
$$
D_j U_n \left( \theta \right) = \frac{2}{J} \int_0^1 D\psi_j(\theta_j,G_{n,j}^{-1})
\bigg\{\psi_j(\theta_j,G_{n,j}^{-1}) - {\frac 1 J  \sum_{k=1}^J} \psi_k(\theta_k,G_{n,k}^{-1})\bigg\},
$$
\begin{equation}\label{eq:derivatives}\tag{A.9}
D_{p,q}  U_n ( \theta ) =- \frac{2}{J^2} \int_0^1 D\psi_p(\theta_p,G_{n,p}^{-1})D\psi_q(\theta_q,G_{n,q}^{-1}), \; p\ne q 
\end{equation}
and
$$
D_{p,p}  U_n (\theta) = \frac{2}{J} \int_0^1 D^2\psi_p(\theta_p,G_{n,p}^{-1})
\bigg\{\psi_j(\theta_j,G_{n,j}^{-1})- \frac 1 J  \sum_{k=1}^J \psi_k ( \theta_k,G_{n,k}^{-1}) \bigg\}
+ \frac{2( J-1) }{J^2} \int_0^1 \left\{D\psi_p(\theta_p,G_{n,p}^{-1}) \right\}^2.  
$$
Using also (\ref{uniformmoments}) we obtain similar expressions for the derivatives of $U(\theta)$, replacing everywhere $G_{n,j}^{-1}$ with $G^{-1}$. We write $DU_n(\theta) = \{ D_j U_n(\theta)\}_{1\leq j\leq J}$, $DU(\theta) = \{D_j U(\theta)\}_{1\leq j\leq J}$ for the gradients and $\Sigma_n(\theta) = \{D_{p,q}U_n ( \theta )\}_{1\leq p,q\leq J}$, $\Sigma(\theta) = \{D_{p,q}U (\theta )\}_{1\leq p,q\leq J}$for the Hessians of $U_n$ and $U$. Note that $\Sigma^*=\Sigma(\theta^*)$ is assumed to be invertible.

We write now $\rho_{n,j}$ for the quantile process based on the $\varepsilon_{i,j}$s.  Observe that \eqref{hyp:tassio13_eps} ensures that we can assume, without loss of generality, that there exist independent Brownian bridges, $B_{n,j}$, satisfying (\ref{csorgohorvath}). Now, recalling that $\psi_j(\theta_j^*,x)=x$, we see that
\begin{equation*}
\sqrt{n} \, D_j U_n( \theta^*)  =  \frac{2}{J} 
\int_0^1  D \psi_j\{\theta_j^*,G_{n,j}^{-1}(t)\}\, \frac{\rho_{n,j}(t)-  \sum_{k=1}^J\rho_{n,k}(t)/J }{g \{G^{-1} (t)\} }dt.\label{eq:expr_D}
\end{equation*}
Now, using (\ref{uniformmoments}) and arguing as in the proof of Theorem~\ref{th:testgen}, we conclude
that 
\begin{align*}
\Bigg|\int_0^{1}  D\psi_j\{\theta_j^*,G_{n,j}^{-1}(t)\}\frac{\rho_{n,k}(t)}{g \left\{ G^{-1} (t) \right\} }
dt -\int_0^{1}  D\psi_j\{\theta_j^*,G^{-1}(t)\} \frac{B_{n,k}(t)}{g \left\{ G^{-1} (t) \right\} } dt\Bigg|  \rightarrow 0
\end{align*}
in probability and, consequently, 
\begin{equation}\label{gradconvergence}\tag{A.10}
\Bigg|\sqrt{n}\, D_j U_n( \theta^*)
-\frac 2 J\int_0^{1}  D\psi_j\{\theta_j^*,G^{-1}(t)\} \frac{B_{n,j}(t) - \sum_{k=1}^J B_{n,k}(t)/J}{g \left\{ G^{-1} (t) \right\} } dt\Bigg|  \rightarrow 0
\end{equation}
in probability.

A further Taylor expansion of $D_j U_n$ around $\theta^*$ shows that for some 
$\tilde{\theta }^n_{j}$ between  $\hat{\theta}_n$ and $\theta^*$ we have
$$
D_jU_n( \hat{\theta }_n)= D_jU_n( \theta ^*) + (D_{1j}U_n (\tilde{\theta }^n_{j}), \ldots, D_{Jj}^2U_n  (\tilde{\theta }^n_{j}) )\times ( \hat{\theta}_n - \theta^* )
$$
and because $\hat{\theta}_n$ is a zero of $DU_n$, we obtain
$$
-D_jU_n( \theta ^*)= (D_{1j}U_n (\tilde{\theta }^n_{j}), \ldots, D_{Jj}U_n  (\tilde{\theta }^n_{j}) )\times ( \hat{\theta}_n - \theta^* ).
$$

Writing $\tilde{\Sigma}_n$ for the $J\times J$ matrix whose  {$J$th}
row equals $(D_{1j}U_n (\tilde{\theta }^n_{j}), \ldots, D_{Jj}U_n 
(\tilde{\theta }^n_{j}) )$, with $j \in \{ 1,\ldots,J\}$, we can rewrite the last expansion as
\begin{equation}
\label{eq:eqbase}\tag{A.11}
- \sqrt{n}\,DU_n( \theta^*) =\tilde{\Sigma}_n \sqrt{n}(\hat{\theta}_n - \theta^*).
\end{equation}
Now, recalling~\eqref{eq:derivatives}, assumptions (\ref{jointsmoothness}) and (\ref{uniformmoments})
yield that $\tilde{\Sigma}_n\to \Sigma^*=\Sigma(\theta^*)$ in probability. As a consequence, 
\eqref{eq:eqbase} and  \eqref{gradconvergence} together with Slutsky's 
Theorem complete the proof of the second claim.

Finally, for the proof of the last claim, since $DU_n(\hat{\theta}_n)=0$, a Taylor expansion around
$\hat{\theta}_n$ shows that 
\begin{equation*}\label{secondorderTaylor}
nU_n(\theta^*)-nU_n(\hat{\theta}_n)=\frac 1 2 \, \{ \sqrt{n} \, (\hat{\theta}_n-\theta^*)\} ^\top \Sigma(\tilde{\theta}_n)
\{\sqrt{n} \, (\hat{\theta}_n-\theta^*)\}
\end{equation*}
for some $\tilde{\theta}_n$ between $\hat{\theta}_n$ and $\theta^*$. Arguing as above we see that 
$\Sigma(\tilde{\theta}_n)\to \Sigma^*$ in probability. Hence, to complete the
proof if suffices to show that 
$$
n U_n(\theta^*)-\frac 1 J \sum_{j=1}^k \int_0^1 \frac {\big\{B_{n,j}(t) - \sum_{k=1}^JB_{n,k}(t) /J \big\}^2}{g\{G^{-1} (t)\}^2}dt \to 0
$$
in probability. Since
$$
n U_n(\theta^*)=\frac 1 J \sum_{j=1}^k \int_0^1 \frac {\big\{\rho_{n,j}(t) - \sum_{k=1}^J\rho_{n,k}(t)/J \big\}^2}{g\{G^{-1} (t)\}^2}dt,
$$
this amounts to proving that 
$$ 
\int_0^1 \ {\{ \rho_{n,j}(t)-B_{n,j}(t) \}^2}/{g \{ G^{-1}(t) \} ^2}dt\to 0
$$
in probability.
Taking $\nu\in(0, 1/ 2)$ in (\ref{csorgohorvath}), we see that 
$$ 
\int_{{1/n} }^{1- {1}/{n} }   {\{\rho_{n,j}(t)-B_{n,j}(t)\}^2} \ {g\{ G^{-1}(t)\}^2}dt\leq O_P(1) \frac 1{n^{1-2\nu}}
\int_{{1/n} }^{1- {1}/{n} }  { \{ t(1-t)\}^{2\nu}}/{g\{G^{-1}(t)\}^2}\to 0,
$$
using condition \eqref{hyp:tassioA3} and dominated convergence. From \eqref{hyp:tassioA3}
we also see that 
$$
\int_{1- {1}/{n} }^1 {B_{n,j}(t)^2}/{g\{G^{-1}(t)\}^2}dt\to 0
$$ 
in probability. Condition~\eqref{hyp:tassioA3} implies also that 
$$
\int_{1- {1}/{n} }^1  {\rho_{n,j}(t)^2}/{g \{G^{-1}(t)\}^2}dt\to 0
$$
in probability; see \cite{SamworthJohnson}. Similar considerations apply to the left tail and complete the proof.
\hfill $\Box$

\bigskip
\noindent
\textbf{Proof of Corollary~\ref{cor:bootstrap_th}.} Writing $\mathcal{L}^*$ for the conditional law given the
$X_{i,j}$s, we see from Theorem~\ref{th:bootstrap} that
$$
W_2^2 \Big[ \mathcal{L} \big[ \sqrt{m_n} \{A_{m_n}(\Theta)\}^{1/2}], \mathcal{L^*} [ \sqrt{m_n} \{A^*_{m_n}(\Theta)\}^{1/2} \big] \Big] \leq
L \frac {m_n}n \frac 1 J \sum_{j=1}^J nW_2^2(\mu,\tilde{\mu}_{n,j}),
$$
where $L=\sup_{\lambda,x,j}\psi_j'(\lambda,x)$, $\mu$ denotes the law of the errors, $\varepsilon_{i,j}$, and 
$\tilde{\mu}_{n,j}$ the empirical distribution function on $\varepsilon_{1,j},\ldots,\varepsilon_{n,j}$. Note that $L<\infty$ by (\ref{hyp:regul_th}),
while $nW_2^2(\mu,\tilde{\mu}_{n,j})=O_P(1)$ as in the proof of Theorem~\ref{prop:tcl}. Hence, we conclude that
\begin{align*}
m_n A^*_{m_n}(\Theta) \rightsquigarrow  
\frac 1 J\sum_{j=1}^J \int_0^1  \left( {\tilde{B}_j}/{g \circ G^{-1}} \right)^2 -  Y^\top \Sigma^{-1}Y/2
\end{align*}
in probability. The conclusion now follows from Lemma 1 in \cite{JanssenPauls}.
\hfill $\Box$

\bigskip
If centering were necessary and we had (\ref{centeredCLT}) rather than the limit in Theorem~\ref{prop:tcl}, we could adapt the last argument as follows. If $A$ and $B$ are positive random variables, then ${\rm E} |A-B|\leq \E (A^{1/2}-B^{1/2})^2+ 2\{ \E A\, \E(A^{1/2}-B^{1/2})^2\}^{1/2}$. We can apply this bound to (an optimal coupling of) $m_nA_{m_n}(\Theta)$ and $m_nA^*_{m_n}(\Theta)$. Now if the errors have a log-concave distribution then $n \E W_2^2(\mu,\tilde{\mu}_{n,j})=O(\ln n)$; see Corollary~6.12 in~\cite{BobkovLedoux}. We conclude that 
$$
W_1[ \mathcal{L} \{ m_n A_{m_n}(\Theta)-c_{m_n}\}, \mathcal{L^*} \{ m_n A^*_{m_n}(\Theta)-c_{m_n}\} ] = W_1 [ \mathcal{L} \{ m_n A_{m_n}(\Theta) \}, \mathcal{L^*} \{ m_n A^*_{m_n}(\Theta)\} ]
$$
vanishes in probability if $m_n=O(n^{\rho})$ for some $\rho \in (0,1)$ . As a consequence,
\begin{align*}
m_n A^*_{m_n}(\Theta) -c_{m_n} \rightsquigarrow  
\frac 1 J\sum_{j=1}^J \int_0^1 \frac{\tilde{B}^2_j-E\tilde{B}^2_j}{(g \circ G^{-1})^2  } -\frac 1 2 Y^{\top}\Sigma^{-1}Y
\end{align*}
in probability.\hfill $\Box$
	
\begin{table}
 \centering
 \caption{Simulations under $\mathcal{H}_0$.\label{tab:frecuenciabajomodelo}} 	
\begin{tabular}{c|ccccccc}
	  		   		$J$&$n$ & $m_{n}=n^{0.6}$ & $m_{n}=n^{0.7}$ & $m_{n}=n^{0.8}$ & $m_{n}=n^{0.9}$ & $m_{n}=n^{0.95}$ & $m_{n}=n$ \\ 
	  		   		\midrule 
	  		   		&50 & 0.144 & 0.079 & 0.038 & 0.046 & 0.041 & 0.03 \\ 
	  		   		\cmidrule{3-8}
	  		   		&100 & 0.148 & 0.067 & 0.07 & 0.05 & 0.04 & 0.033 \\ 
	  		   		\cmidrule{3-8}
	  		   		&200 & 0.129 & 0.085 & 0.068 & 0.043 & 0.037 & 0.044 \\ 
	  		   		\cmidrule{3-8}
	  		   		2&500 & 0.138 & 0.089 & 0.05 & 0.048 & 0.035 & 0.036 \\ 
	  		   		\cmidrule{3-8}
	  		   		&1000 & 0.127 & 0.086 & 0.063 & 0.055 & 0.039 & 0.032 \\
	  		   		\cmidrule{3-8} 
	  		   		&2000 & 0.129 & 0.104 & 0.071 & 0.048 & 0.043 & 0.038 \\
	  		   		\cmidrule{3-8} 
	  		   		&5000 & 0.039 & 0.042 & 0.041 & 0.049 & 0.043 & 0.055\\
	  		   		\midrule 
	  		   
	  		   		&50 & 0.295 & 0.194 & 0.115 & 0.078 & 0.054 & 0.034 \\
	  		   		\cmidrule{3-8} 
	  		   		&100 & 0.273 & 0.163 & 0.089 & 0.053 & 0.034 & 0.039 \\ 
	  		   		\cmidrule{3-8}
	  		   		&200 & 0.238 & 0.15 & 0.077 & 0.054 & 0.047 & 0.031 \\ 
	  		   		\cmidrule{3-8}
	  		   		3&500 & 0.226 & 0.122 & 0.07 & 0.057 & 0.042 & 0.029 \\ 
	  		   		\cmidrule{3-8}
	  		   		&1000 & 0.217 & 0.107 & 0.092 & 0.069 & 0.042 & 0.035 \\ 
	  		   		\cmidrule{3-8}
	  		   		&2000 & 0.221 & 0.128 & 0.077 & 0.053 & 0.043 & 0.035 \\ 
	  		   		\cmidrule{3-8}
	  		   		&5000 & 0.205 & 0.145 & 0.082 & 0.06 & 0.025 & 0.047 \\
	  		   		\midrule 
	  		   		  		   	
	  		   		&50 & 0.659 & 0.428 & 0.281 & 0.129 & 0.111 & 0.081 \\
	  		   		\cmidrule{3-8} 
	  		   		&100 & 0.583 & 0.337 & 0.192 & 0.104 & 0.083 & 0.053 \\ 
	  		   		\cmidrule{3-8}
	  		   		&200 & 0.538 & 0.281 & 0.159 & 0.081 & 0.078 & 0.029 \\ 
	  		   		\cmidrule{3-8}
	  		   		5&500 & 0.449 & 0.267 & 0.138 & 0.063 & 0.056 & 0.04 \\ 
	  		   		\cmidrule{3-8}
	  		   		&1000 & 0.415 & 0.238 & 0.129 & 0.064 & 0.051 & 0.037 \\
	  		   		\cmidrule{3-8} 
	  		   		&2000 & 0.354 & 0.212 & 0.115 & 0.06 & 0.053 & 0.032 \\ 
	  		   		\cmidrule{3-8}
	  		   		&5000 & 0.322 & 0.203 & 0.108 & 0.057 & 0.061 & 0.039 \\
	  		   		\midrule 
	  		   		
	  		   		&50 & 0.996 & 0.971 & 0.873 & 0.702 & 0.553 & 0.456 \\ 
	  		   		\cmidrule{3-8}
	  		   		&100 & 0.994 & 0.902 & 0.708 & 0.433 & 0.33 & 0.226 \\
	  		   		\cmidrule{3-8} 
	  		   		&200 & 0.958 & 0.802 & 0.521 & 0.247 & 0.184 & 0.119 \\
	  		   		\cmidrule{3-8} 
	  		   		10&500 & 0.914 & 0.663 & 0.388 & 0.149 & 0.093 & 0.063 \\
	  		   		\cmidrule{3-8} 
	  		   		&1000 & 0.864 & 0.532 & 0.286 & 0.119 & 0.084 & 0.046 \\ 
	  		   		\cmidrule{3-8}
	  		   		&2000 & 0.813 & 0.473 & 0.239 & 0.103 & 0.063 & 0.051 \\ 
	  		   		\cmidrule{3-8}
	  		   		&5000 & 0.756 & 0.449 & 0.217 & 0.088 & 0.061 & 0.041 \\
	  		   		\bottomrule
	  		   		\end{tabular}
	  		   	\end{table}  		   		  	
	  		   	
	  		   	\begin{table}
\caption{Power of the test for $\gamma {=}_d \varepsilon (1 )$.\label{tab:potenciaexp}}		
	  		   		\centering   	
	  		   		\begin{tabular}{c|ccccccc}
	  		   		$J$&$n$ & $m_{n}=n^{0.6}$ & $m_{n}=n^{0.7}$ & $m_{n}=n^{0.8}$ & $m_{n}=n^{0.9}$ & $m_{n}=n^{0.95}$& $m_{n}=n$ \\
	  		   		\toprule 
	  		   		&50 & 0.961 & 0.919 & 0.897 & 0.864 & 0.829 & 0. 838\\ 
	  		   		\cmidrule{3-8}
	  		   		&100 & 1 & 0.998 & 0.998 & 0.995 & 0.994 & 0.993\\ 
	  		   		\cmidrule{3-8}
	  		   		&200 & 1 & 1 & 1 & 1 & 1 &1 \\ 
	  		   		\cmidrule{3-8}
	  		   		2&500 & 1 & 1 & 1 & 1 & 1 & 1 \\ 
	  		   		\cmidrule{3-8}
	  		   		&1000 & 1 & 1 & 1 & 1 & 1 & 1 \\ 
	  		   		\cmidrule{3-8}
	  		   		&2000 & 1 & 1 & 1 & 1 & 1 & 1 \\ 
	  		   		\cmidrule{3-8}
	  		   		&5000 & 1 & 1 & 1 & 1 & 1 & 1\\
	  		   		\midrule 
	  		   		&50 & 0.987 & 0.971 & 0.97 & 0.953 & 0.939 & 0.91 \\ 
	  		   		\cmidrule{3-8}
	  		   		&100 & 1 & 1 & 0.999 & 1 & 0.999 & 0.999 \\ 
	  		   		\cmidrule{3-8}
	  		   		&200 & 1 & 1 & 1 & 1 & 1 & 1\\ 
	  		   		\cmidrule{3-8}
	  		   		3&500 & 1 & 1 & 1 & 1 & 1 & 1 \\ 
	  		   		\cmidrule{3-8}
	  		   		&1000 & 1 & 1 & 1 & 1 & 1 & 1 \\ 
	  		   		\cmidrule{3-8}
	  		   		&2000 & 1 & 1 & 1 & 1 & 1 & 1 \\
	  		   		\cmidrule{3-8} 
	  		   		&5000 & 1 & 1 & 1 & 1 & 1 & 1 \\
	  		   		\midrule 
	  		   		&50 & 1 & 0.996 & 0.988 & 0.976 & 0.971 & 0.955 \\ 
	  		   		\cmidrule{3-8}
	  		   		&100 & 1 & 1 & 1 & 1 & 1 & 1 \\
	  		   		\cmidrule{3-8}
	  		   		&200 & 1 & 1 & 1 & 1 & 1 & 1\\ 
	  		   		\cmidrule{3-8}
	  		   		5&500 & 1 & 1 & 1 & 1 & 1 & 1 \\ 
	  		   		\cmidrule{3-8}
	  		   		&1000 & 1 & 1 & 1 & 1 & 1 & 1 \\
	  		   		\cmidrule{3-8} 
	  		   		&2000 & 1 & 1 & 1 & 1 & 1 & 1 \\
	  		   		\cmidrule{3-8} 
	  		   		&5000 & 1 & 1 & 1 & 1 & 1 & 1\\
	  		   		\midrule 
	  		   		&50 & 1 & 1 & 1 & 1 & 0.996 & 0.985 \\
	  		   		\cmidrule{3-8}
	  		   		&100 & 1 & 1 & 1 & 1 & 1 & 1 \\ 
	  		   		\cmidrule{3-8}
	  		   		&200 & 1 & 1 & 1 & 1 & 1 & 1\\ 
	  		   		\cmidrule{3-8}
	  		   		10&500 & 1 & 1 & 1 & 1 & 1 & 1 \\ 
	  		   		\cmidrule{3-8}
	  		   		&1000 & 1 & 1 & 1 & 1 & 1 & 1 \\ 
	  		   		\cmidrule{3-8}
	  		   		&2000 & 1 & 1 & 1 & 1 & 1 & 1 \\ 
	  		   		\cmidrule{3-8}
	  		   		&5000 & 1 & 1 & 1 & 1 & 1 & 1\\
	  		   		\bottomrule
	  		   	\end{tabular}
	  		   	\end{table}

	  		   	\begin{table}
					  		   \caption{Power  of the test $\gamma \stackrel{d}{=} \mbox{Laplace}\left(0,1\right).$
			    \label{tab:potenciadobleexp}}
	  		   		\centering
	  		   	\begin{tabular}{c|ccccccc}
	  		   	\toprule
	  		   	$J$& $n$ & $m_{n}=n^{0.6}$ & $m_{n}=n^{0.7}$ & $m_{n}=n^{0.8}$ & $m_{n}=n^{0.9}$ & $m_{n}=n^{0.95}$ & $m_{n}=n$  \\
	  		   	\midrule 
	  		   	&50  & 0.426 & 0.33 & 0.3 & 0.241 & 0.223 & 0.163\\
	  		   	\cmidrule{3-8}
	  		   	&100 & 0.658 & 0.534 & 0.468 & 0.365 & 0.361 & 0.3 \\ 
	  		   	\cmidrule{3-8}
	  		   	&200 & 0.855 & 0.824 & 0.751 & 0.665 & 0.613 & 0.602\\
	  		   	\cmidrule{3-8} 
	  		   	2&500 & 0.998 & 0.998 & 0.993 & 0.982 & 0.965 &0.962 \\ 
	  		   	\cmidrule{3-8}
	  		   	&1000 & 1 & 1 & 1 & 1 & 0.999 & 1 \\ 
	  		   	\cmidrule{3-8}
	  		   	&2000 & 1 & 1 & 1 & 1 & 1 & 1 \\ 
	  		   	\cmidrule{3-8}
	  		   	&5000 & 1 & 1 & 1 & 1 & 1 & 1 \\
	  		   	\midrule 
	  		   	&50 & 0.657 & 0.533 & 0.422 & 0.331 & 0.282 & 0.223 \\
	  		   	\cmidrule{3-8}
	  		   	&100 & 0.831 & 0.708 & 0.586 & 0.514 & 0.461 & 0.377 \\ 
	  		   	\cmidrule{3-8}
	  		   	&200  & 0.946 & 0.915 & 0.841 & 0.778 & 0.709 & 0.661\\ 
	  		   	\cmidrule{3-8}
	  		   	3&500 & 1 & 0.998 & 0.997 & 0.994 & 0.989 & 0.977 \\ 
	  		   	\cmidrule{3-8}
	  		   	&1000 & 1 & 1 & 1 & 1 & 1 & 1 \\ 
	  		   	\cmidrule{3-8}
	  		   	&2000 & 1 & 1 & 1 & 1 & 1 & 1 \\ 
	  		   	\cmidrule{3-8}
	  		   	&5000 & 1 & 1 & 1 & 1 & 1 & 1 \\
	  		   	\midrule 
	  		   	&50 & 0.895 & 0.741 & 0.633 & 0.471 & 0.394 & 0.333 \\ 
	  		   	\cmidrule{3-8}
	  		   	&100 & 0.936 & 0.874 & 0.728 & 0.623 & 0.519 & 0.443\\ 
	  		   	\cmidrule{3-8}
	  		   	&200 & 0.994 & 0.947 & 0.903 & 0.847 & 0.786 & 0.696\\ 
	  		   	\cmidrule{3-8}
	  		   	5&500 & 1 & 1 & 1 & 0.996 & 0.992 & 0.985 \\ 
	  		   	\cmidrule{3-8}
	  		    &1000 & 1 & 1 & 1 & 1 & 1 & 1 \\ 
	  		    \cmidrule{3-8}
	  		   	&2000 & 1 & 1 & 1 & 1 & 1 & 1 \\ 
	  		   	\cmidrule{3-8}
	  		   	&5000 & 1 & 1 & 1 & 1 & 1 & 1\\
	  		   	\midrule 
	  		   	&50 & 1 & 0.997 & 0.97 & 0.875 & 0.79 & 0.703 \\
	  		   	\cmidrule{3-8}
	  		   	&100 & 0.997 & 0.985 & 0.949 & 0.854 & 0.765 & 0.643\\
	  		   	\cmidrule{3-8} 
	  		   	&200 & 1 & 0.996 & 0.968 & 0.924 & 0.859 & 0.789\\ 
	  		   	\cmidrule{3-8}
	  		   	10&	500 & 1 & 1 & 1 & 0.996 & 0.996 & 0.975 \\ 
	  		   	\cmidrule{3-8}
	  		   	&1000 & 1 & 1 & 1 & 1 & 1 & 0.999\\ 
	  		   	\cmidrule{3-8}
	  		   	&2000 & 1 & 1 & 1 & 1 & 1 & 1 \\ 
	  		   	\cmidrule{3-8}
	  		   	&5000 & 1 & 1 & 1 & 1 & 1 & 1\\
	  		   	\bottomrule
	  		   \end{tabular}

	  		   \end{table}

	  		   	    \begin{table}
				    	  		   		\caption{Power of the test $\gamma \stackrel{d}{=} t_{(3)}.$\label{tab:potenciat3}}

	  		   	    	\centering	
	  		   		\begin{tabular}{c|ccccccc}
	  		   		\toprule
	  		   		$I$&$n$ & $m_{n}=n^{0.6}$ & $m_{n}=n^{0.7}$ & $m_{n}=n^{0.8}$ & $m_{n}=n^{0.9}$ & $m_{n}=n^{0.95}$ & $m_{n}=n$\\
	  		   		\midrule 
	  		   		&50 & 0.566 & 0.445 & 0.429 & 0.352 & 0.321 & 0.307\\
	  		   		\cmidrule{3-8}
	  		   		&100 & 0.775 & 0.704 & 0.647 & 0.576 & 0.503 & 0.454\\
	  		   		\cmidrule{3-8}
	  		   		&200  & 0.942 & 0.927 & 0.882 & 0.833 & 0.771 & 0.697 \\ 
	  		   		\cmidrule{3-8}
	  		   		2&500 & 1 & 0.997 & 0.995 & 0.991 & 0.989 & 0.957 \\ 
	  		   		\cmidrule{3-8}
	  		   		&1000 & 1 & 1 & 1 & 1 & 1 & 0.986 \\ 
	  		   		\cmidrule{3-8}
	  		   		&2000 & 1 & 1 & 1 & 1 & 1 &0.999 \\ 
	  		   		\cmidrule{3-8}
	  		   		&5000 & 1 & 1 & 1 & 1 & 1 & 0.997\\
	  		   		\midrule 	  		   		
	  		   		&50 & 0.745 & 0.653 & 0.546 & 0.46 & 0.402 & 0.349 \\ 
	  		   		\cmidrule{3-8}
	  		   		&100 & 0.881 & 0.821 & 0.738 & 0.65 & 0.592 & 0.563 \\
	  		   		\cmidrule{3-8} 
	  		   		&200 & 0.98 & 0.958 & 0.928 & 0.891 & 0.873 & 0.794 \\
	  		   		\cmidrule{3-8}
	  		   		3&500 & 1 & 1 & 0.999 & 0.997 & 0.997 & 0.978 \\ 
	  		   		\cmidrule{3-8}
	  		   		&1000 & 1 & 1 & 1 & 1 & 1 & 0.995\\ 
	  		   		\cmidrule{3-8}
	  		   		&2000 & 1 & 1 & 1 & 1 & 1 & 1 \\ 
	  		   		\cmidrule{3-8}
	  		   		&5000 & 1 & 1 & 1 & 1 & 1 & 1\\
	  		   		\midrule   		   	  		   		
	  		   		&50  & 0.91 & 0.813 & 0.682 & 0.593 & 0.525 & 0.45\\ 
	  		   		\cmidrule{3-8}
	  		   		&100 & 0.972 & 0.909 & 0.822 & 0.751 & 0.686 & 0.621\\
	  		   		\cmidrule{3-8}
	  		   		&200 & 0.995 & 0.984 & 0.967 & 0.915 & 0.887 & 0.836 \\
	  		   		\cmidrule{3-8} 
	  		   		5&500 & 1 & 1 & 1 & 0.999 & 0.999 & 0.995\\ 
	  		   		\cmidrule{3-8}
	  		   		&1000 & 1 & 1 & 1 & 1 & 1 & 1 \\
	  		   		\cmidrule{3-8}
	  		   		&2000 & 1 & 1 & 1 & 1 & 1 & 1 \\ 
	  		   		\cmidrule{3-8}
	  		   		&5000 & 1 & 1 & 1 & 1 & 1 & 1\\
	  		   		\midrule 
	  		   		&50 & 1 & 0.997 & 0.953 & 0.894 & 0.827 & 0.758 \\
	  		   		\cmidrule{3-8}
	  		   		&100 & 0.999 & 0.993 & 0.969 & 0.907 & 0.862 & 0.79\\ 
	  		   		\cmidrule{3-8}
	  		   		&200 & 1 & 0.998 & 0.995 & 0.961 & 0.941 & 0.903 \\ 
	  		   		\cmidrule{3-8}
	  		   		10&500 & 1 & 1 & 1 & 1 & 0.998 & 0.988\\ 
	  		   		\cmidrule{3-8}
	  		   		&1000 & 1 & 1 & 1 & 1 & 1 & 0.998 \\ 
	  		   		\cmidrule{3-8}
	  		   		&2000 & 1 & 1 & 1 & 1 & 1 & 0.999 \\ 
	  		   		\cmidrule{3-8}
	  		   		&5000 & 1 & 1 & 1 & 1 & 1 & 1 \\
	  		   		\bottomrule
	  		   		\end{tabular}
	  		   		\end{table}
	  		   		
	  		   		\begin{table}
						  		   		\caption{Power of the test $ \gamma \stackrel{d}{=} t_{(4)}.$\label{tab:potenciat4}}

	  		   			\centering
	  		   		\begin{tabular}{c|ccccccc}
	  		   		\toprule
	  		   		I&$n$ & $m_{n}=n^{0.6}$ & $m_{n}=n^{0.7}$ & $m_{n}=n^{0.8}$ & $m_{n}=n^{0.9}$ & $m_{n}=n^{0.95}$ & $m_{n}=n$ \\
	  		   		\midrule 
	  		   		&50 & 0.398 & 0.353 & 0.292 & 0.207 & 0.182 & 0.183 \\
	  		   		\cmidrule{3-8}
	  		   		&100 & 0.623 & 0.52 & 0.429 & 0.341 & 0.29 & 0.228\\ 
	  		   		\cmidrule{3-8}
	  		   		&200 & 0.826 & 0.717 & 0.65 & 0.589 & 0.526 & 0.41 \\
	  		   		\cmidrule{3-8}
	  		   		2&500 & 0.989 & 0.978 & 0.954 & 0.928 & 0.878 & 0.787 \\
	  		   		\cmidrule{3-8} 
	  		   		&1000 & 1 & 1 & 0.999 & 1 & 0.984 & 0.955\\ 
	  		   		\cmidrule{3-8}
	  		   		&2000 & 1 & 1 & 1 & 1 & 1 & 0.985\\ 
	  		   		\cmidrule{3-8}
	  		   		&5000 & 1 & 1 & 1 & 1 & 1 & 0.993\\
	  		   		\midrule 
	  		   		&50 & 0.634 & 0.495 & 0.4 & 0.295 & 0.263 & 0.222 \\
	  		   		\cmidrule{3-8}
	  		   		&100 & 0.756 & 0.666 & 0.56 & 0.465 & 0.399 & 0.336\\ 
	  		   		\cmidrule{3-8}
	  		   		&200  & 0.914 & 0.859 & 0.778 & 0.663 & 0.602 & 0.521 \\
	  		   		\cmidrule{3-8}
	  		   		3&500 & 0.998 & 0.989 & 0.985 & 0.972 & 0.928 & 0.868\\ 
	  		   		\cmidrule{3-8}
	  		   		&1000 & 1 & 1 & 1 & 1 & 0.999 & 0.963 \\
	  		   		\cmidrule{3-8}
	  		   		&2000 & 1 & 1 & 1 & 1 & 1 & 0.989 \\ 
	  		   		\cmidrule{3-8}
	  		   		&5000 & 1 & 1 & 1 & 1 & 1 & 1\\
	  		   		\midrule 
	  		   		&50  & 0.851 & 0.709 & 0.583 & 0.426 & 0.359 & 0.316\\ 
	  		   		\cmidrule{3-8}
	  		   		&100 & 0.919 & 0.825 & 0.668 & 0.546 & 0.493 & 0.316 \\
	  		   		\cmidrule{3-8}
	  		   		&200 & 0.959 & 0.908 & 0.842 & 0.738 & 0.684 & 0.578 \\ 
	  		   		\cmidrule{3-8}
	  		   		5&500 & 1 & 0.997 & 0.994 & 0.973 & 0.934 & 0.888 \\ 
	  		   		\cmidrule{3-8}
	  		   		&1000 & 1 & 1 & 1 & 1 & 0.999 & 0.968\\ 
	  		   		\cmidrule{3-8}
	  		   		&2000 & 1 & 1 & 1 & 1 & 1 & 1 \\ 
	  		   		\cmidrule{3-8}
	  		   		&5000 & 1 & 1 & 1 & 1 & 1 & 0.999\\
	  		   		\midrule 
	  		   		&50 & 1 & 0.986 & 0.941 & 0.813 & 0.774 & 0.653\\ 
	  		   		\cmidrule{3-8}
	  		   		&100 & 1 & 0.988 & 0.925 & 0.806 & 0.738 & 0.606 \\ 
	  		   		\cmidrule{3-8}
	  		   		&200 & 1 & 0.991 & 0.948 & 0.854 & 0.813 & 0.679 \\
	  		   		\cmidrule{3-8} 
	  		   		10&500 & 1 & 1 & 0.998 & 0.985 & 0.954 & 0.886\\ 
	  		   		\cmidrule{3-8}
	  		   		&1000 & 1 & 1 & 1 & 1 & 0.997 & 0.949 \\ 
	  		   		\cmidrule{3-8}
	  		   		&2000 & 1 & 1 & 1 & 1 & 1 & 0.974\\ 
	  		   		\cmidrule{3-8}
	  		   		&5000 & 1 & 1 & 1 & 1 & 1 & 0.995\\
	  		   		\bottomrule
	  		   		\end{tabular} 
	  		   		\end{table}
{ 
\section*{Ackowledgements}
The authors would like to thank the Editor for his careful reading of the manuscript and his suggestions that improved the paper.

E. del Barrio's research has been partially supported by Spanish Ministerio 
de Econom\'ia y Competitividad, grant MTM2017-86061-C2-1-P, and by Consejer\'ia de 
Educaci\'on de la Junta de Castilla y Le\'on and FEDER, grants VA005P17 and 
VA002G18}

\end{document}